\documentclass[a4paper, 11pt]{article}

\usepackage{inputenc}
\usepackage{amsmath, amsthm, amssymb}
\usepackage{color}
\usepackage{a4wide}
\usepackage{graphicx}

\numberwithin{equation}{section}

\newtheorem{thm}{Theorem}[section]

\newtheorem{cor}[thm]{Corollary}
\newtheorem{defi}[thm]{Definition}

\theoremstyle{definition}

\newtheorem{rem}[thm]{Remark}

\newenvironment{prf}{\noindent\textbf{Proof.}}{\hfill$\square$\\[-0.5ex]}

\newenvironment{abst}{\begin{minipage}[c]{0.9\textwidth} \footnotesize \textbf{Abstract.}}
{\end{minipage}\\[2ex]}

\newenvironment{key}{\begin{minipage}[c]{0.9\textwidth} \footnotesize \textbf{Keywords.}}
{\end{minipage}\\[2ex]}

\newcommand{\eps}{\ensuremath \varepsilon}
 
\DeclareMathOperator{\dd}{\ensuremath\normalfont{d}}
\DeclareMathOperator{\DD}{\ensuremath\normalfont{D}}
\DeclareMathOperator{\CC}{\ensuremath\normalfont{C}}
\DeclareMathOperator{\LL}{\ensuremath\normalfont{L}}

\newcommand{\g}{G^{(2)}}
\newcommand{\ooox}{{\bar{o}}^{(1)}}
\newcommand{\oooy}{{\bar{o}}^{(2)}}
\newcommand{\oooz}{{\bar{o}}^{(3)}}
\newcommand{\oox}{\tilde{o}^{(1)}}
\newcommand{\ooy}{\tilde{o}^{(2)}}
\newcommand{\ooz}{\tilde{o}^{(3)}}

\newcommand{\ox}{{o}^{(1)}}
\newcommand{\oy}{{o}^{(2)}}
\newcommand{\oz}{{o}^{(3)}}

\newcommand{\Fx}{\tilde{F}_1}
\newcommand{\Fy}{\tilde{F}_2}
\newcommand{\Fz}{\tilde{F}_3}
\newcommand{\fx}{\tilde{f}_1}
\newcommand{\fy}{\tilde{f}_2}
\newcommand{\fz}{\tilde{f}_3}
\newcommand{\Ffx}{{\bar{F}}_1}
\newcommand{\Ffy}{{\bar{F}}_2}
\newcommand{\Ffz}{{\bar{F}}_3}
\newcommand{\ffx}{{\bar{f}}_1}
\newcommand{\ffy}{{\bar{f}}_2}
\newcommand{\ffz}{{\bar{f}}_3}

\begin{document}
\begin{center}
\Large\bfseries On the Unique Reconstruction of Induced Spherical Magnetizations \normalsize\mdseries
\\[3ex]
{Christian Gerhards}\footnote{University of Vienna, Computational Science Center\\Oskar-Morgenstern-Platz 1, 1090 Vienna\\e-mail: christian.gerhards@univie.ac.at}
\\[3ex]
\today
\end{center}

\begin{abst}
Recovering spherical magnetizations $m$ from magnetic field data in the exterior is a highly non-unique problem. A spherical Hardy-Hodge decomposition supplies information on what contributions of the magnetization $m$ are recoverable but it does not supply geophysically suitable constraints on $m$ that would guarantee uniqueness for the entire magnetization. In this paper, we focus on the case of induced spherical magnetizations and show that uniqueness is guaranteed if one assumes that the magnetization is compactly  supported on the sphere. The results are based on ideas presented in \cite{baratchart13} for the planar setting.
\end{abst}

\begin{key}
Spherical Hardy-Hodge decomposition, potential theory, spherical magnetization
\end{key}


\section{Introduction}

The lithospheric contribution to the Earth's magnetic field is due to magnetized rocks in the Earth's crust, which can be expressed as a spherical shell $\Omega_{r,R}=\{x\in\mathbb{R}^3:r<|x|<R\}$. The magnetic potential $V$ that is generated by a vectorial magnetization $m:\Omega_{r,R}\to\mathbb{R}^3$ can be expressed by
\begin{align}\label{eqn:smag1}
 V(x)=\frac{1}{4\pi}\int_{\Omega_{r,R}}m(y)\cdot\frac{x-y}{|x-y|^3}{\dd}y,\quad x\in\mathbb{R}^3.
\end{align}
Recovering $m$ from knowledge of $V$ only in the exterior $\Omega_R^{ext}=\{x\in\mathbb{R}^3:|x|>R\}$ is a highly non-unique problem. The relation \eqref{eqn:smag1} actually seems to be a vectorial version of the gravimetry problem (see, e.g., \cite{ballani82, ballani93, michel05, michel08}) and reveals similar uniqueness issues. However, opposed to the gravimetry problem, where the assumption of a harmonic mass density leads to uniqueness, the assumption of a harmonic magnetization $m$ would still maintain a certain non-uniqueness. The non-uniqueness even persists if we restrict ourselves to induced magnetizations of the form $m=Qv$, where $v$ denotes a known inducing vector field and $Q$ an unknown scalar susceptibility. In \cite{runcorn75}, it has been shown that a constant susceptibility in the spherical shell $\Omega_{r,R}$ produces no magnetic effect in the exterior $\Omega_R^{ext}$ (for any inducing vector field of the form $v=\nabla U$, where $U$ is harmonic in $\Omega_R^{ext}$). These considerations have been generalized to ellipsoidal shells in \cite{jackson99}. A discussion of further examples of magnetizations and related uniqueness issues can be found, e.g., in \cite{blakely95}.

Since the thickness of the spherical shell $\Omega_{r,R}$, where magnetization in the Earth's lithosphere occurs, is only a few tens of kilometers (thus, negligibly small compared to the Earth's radius), it is geophysically reasonable to reduce the considerations to vertically integrated magnetizations $m:\Omega_R\to\mathbb{R}^3$ on the sphere $\Omega_R=\{x\in\mathbb{R}^3:|x|=R\}$. Therefore, from now on, we consider the relation 
\begin{align}\label{eqn:smag10}
 V(x)=\frac{1}{4\pi}\int_{\Omega_R}m(y)\cdot\frac{x-y}{|x-y|^3}{\dd}\omega(y),\quad x\in\mathbb{R}^3,
\end{align}
instead of \eqref{eqn:smag1}. By ${\dd}\omega$ we denote the surface element on the sphere $\Omega_R$. The non-uniqueness of recovering a vertically integrated magnetization $m$ from the knowledge of $V$ in $\Omega_R^{ext}$ can be characterized by a fairly well-known decomposition (see, e.g., \cite{backus96, freedengerhards12, gerhards11a, gerhards12, gubbins11, mayer06, mayermaier06, olsen10b})
\begin{align}\label{eqn:hardhodge0}
m=\tilde{m}^{(1)}+\tilde{m}^{(2)}+\tilde{m}^{(3)},
\end{align}
which has the property that $V\equiv0$ in $\Omega_R^{ext}$ if and only if $\tilde{m}^{(2)}\equiv0$ (in other words, any magnetization of the form $m=\tilde{m}^{(1)}+\tilde{m}^{(3)}$ produces no magnetic potential in the exterior $\Omega_R^{ext}$). We call such a decomposition a Hardy-Hodge decomposition (cf. \cite{baratchart13} for its Euclidean counterpart in $\mathbb{R}^2$) and treat it in more detail later on. An illustration of this decomposition for recent magnetization models is supplied in \cite{gubbins11}. Nonetheless, a characterization of $m$ by \eqref{eqn:hardhodge0} still states that the contributions $\tilde{m}^{(1)}$ and $\tilde{m}^{(3)}$ cannot be reconstructed from knowledge of $V$ in $\Omega_R^{ext}$. Even if we assume an induced magnetization $m=Qv$, non-trivial susceptibilities $Q$ have been constructed in \cite{maushaak03} that generate a magnetic potential $V$ via \eqref{eqn:smag10} which vanishes in $\Omega_R^{ext}$ (again, assuming that $v$ is known and of the form $v=\nabla U$, where $U$ is harmonic in $\Omega_R^{ext}$), underlining the non-uniqueness of the problem.

In this paper, we show that an induced magnetization $m=Qv$ is uniquely recoverable from \eqref{eqn:smag10} if $V$ is known in $\Omega_R^{ext}$ and if one imposes the additional condition that $m$ has compact support in $\Omega_R$ (it is not necessary that $v$ is of the form $v=\nabla U$). Now, if there exists a model of the vertically integrated induced magnetization that is very accurate in some region $\Gamma_R\subset\Omega_R$ of the Earth, the residuum of the magnetic potential generated by this model magnetization and the magnetic potential obtained from actual global magnetic field measurements forms a magnetic potential that can be regarded as being generated by a magnetization with compact support in $\Omega_R\setminus{\Gamma_R}$. The residual magnetization can then be determined uniquely due to our result and, together with the accurate model magnetization in $\Gamma_R$, we obtain a trustworthy model for the induced vertically integrated magnetization in $\Omega_R$. The proof of our uniqueness result is based on a combination of the spherical Hardy-Hodge decomposition from  \cite{freedengerhards12, gerhards11a, gerhards12} and the ideas presented in \cite{baratchart13} for the case of thin-plate magnetizations in the plane $\mathbb{R}^2$.

After introducing some notations in Section \ref{sec:notation}, we briefly recapitulate results from \cite{baratchart13, lima13} in Section \ref{sec:tpc} in order to highlight their relations to the spherical case later on. In Section \ref{sec:sdecomp}, we introduce the spherical Hardy-Hodge decomposition and in Section \ref{sec:sunique}, we discuss constraints on induced spherical magnetizations $m=Qv$ that guarantee uniqueness if $V$ is known, e.g., only in the exterior $\Omega_R^{ext}$. Furthermore, we discuss the numerical reconstruction of $m$ and supply some examples in Section \ref{sec:num}.

\subsection{Notations}\label{sec:notation}

For the spherical setting, we assume that the magnetization $m$ is located on the unit sphere $\Omega=\Omega_1$ and we denote the exterior by $\Omega^{ext}=\{x\in\mathbb{R}^3:|x|>1\}$ and the interior by $\Omega^{int}=\{x\in\mathbb{R}^3:|x|<1\}$. In the Euclidean setting of thin-plate magnetizations in \cite{baratchart13}, the magnetization $m$ is restricted to $\mathbb{R} ^2\simeq\mathbb{R}^2\times\{0\}$ and the \lq\lq exterior'' is represented by the upper half-space $\mathbb{R}^3_+$ and the \lq\lq interior'' by the lower half-space $\mathbb{R}_-^3$. In the following, we essentially use identical notations for the spherical setting and for the  Euclidean thin-plate setting (of which the thin-plate setting only appears in Section \ref{sec:tpc}). Vector fields mapping $\Omega$ or $\mathbb{R}^2$ into $\mathbb{R}^3$ are denoted by lower case letters $f,g,\ldots$, and if they are square-integrable, we say that they are of class ${\rm l}^2(\Omega)$ or ${\rm l}^2(\mathbb{R}^2)$, respectively. Scalar fields mapping  $\Omega$ or $\mathbb{R}^2$ into $\mathbb{R}$ are denoted by upper case letters $F,G,\ldots$, and if they are square-integrable, we say that they are of class $\LL^2(\Omega)$ or $\LL^2(\mathbb{R}^2)$, respectively. By Latin letters $x,y,\ldots$ we mean vectors in $\mathbb{R}^3$ or $\mathbb{R} ^2$, by Greek letters $\xi,\eta,\ldots$ unit vectors in $\Omega$.

The surface gradient $\nabla^*$ is defined as $\nabla^*_x=\big(\frac{\partial}{\partial{x_1}},\frac{\partial}{\partial{x_2}},0\big)$, for $x=(x_1,x_2)\in\mathbb{R}^2\simeq\mathbb{R}^2\times\{0\}$, in case of the Euclidean setting. In the spherical setting, it is defined by the connection $\nabla^*_\xi=r\nabla_x-r\xi\frac{\partial}{\partial r}$ to the gradient $\nabla_x$ in $\mathbb{R}^3$, where $r=|x|$ and $\xi=\frac{x}{|x|}\in\Omega$. The surface curl gradient $\LL^*$ reads as $\LL^*_x=(0,0,1)\wedge \nabla^*_x=\big(-\frac{\partial}{\partial{x_2}},\frac{\partial}{\partial{x_1}},0\big)$ in the Euclidean setting and as $\LL^*_\xi=\xi\wedge \nabla^*_\xi$ in the spherical setting, where $\wedge$ denotes the vector product. For convenience, we introduce the following notation for the (tangential) differential operators above:
\begin{align}
 o^{(2)}=\nabla^*,\qquad o^{(3)}=\LL^*.
\end{align}
These operators are complemented by the operator $o^{(1)}$  which is given by 
\begin{align}
 o^{(1)}=(0,0,1)\,\textnormal{id}
\end{align}
in the Euclidean setting and by 
\begin{align}
 o^{(1)}_\xi=\xi\,\textnormal{id}
\end{align}
for the spherical setting (id denotes the identity operator). $o^{(1)}$ always points in normal direction with respect to $\mathbb{R}^2\simeq\mathbb{R}^2\times\{0\}$ or $\Omega$, respectively.
Last, we need the Beltrami operator $\Delta^*$ which is defined by $\Delta^*=\nabla^*\cdot \nabla^*=\LL^*\cdot\LL^*$. In the Euclidean setting this means that $\Delta^*_x=\big(\frac{\partial}{\partial x_1}\big)^2+\big(\frac{\partial}{\partial x_2}\big)^2$, for $x=(x_1,x_2)\in\mathbb{R}^2\simeq\mathbb{R}^2\times\{0\}$, and in the spherical setting that the connection $\Delta^*_\xi=r^2\Delta_x-\frac{\partial}{\partial r}r^2 \frac{\partial}{\partial r}$ to the Laplace operator $\Delta_x$ in $\mathbb{R}^3$ holds true, where $r=|x|$ and $\xi=\frac{x}{|x|}\in\Omega$. 


\subsection{The Thin-Plate Case}\label{sec:tpc}
For the thin-plate setting, we assume that the magnetic potential $V$ is generated by a vectorial magnetization $m$ of class ${\rm l}^2(\mathbb{R}^2)$. Then we can write
\begin{align}\label{eqn:mag2d}
 V(x)=\frac{1}{4\pi}\int_{\mathbb{R}^2}m(y)\cdot\frac{x-(y,0)}{|x-(y,0)|^3}{\dd} y,\qquad x\in\mathbb{R}^3.
\end{align}
We use the following definition in order to characterize $m$ with respect to its effect on $V$. 

\begin{defi}\label{def:equivabove}
Two magnetizations $m$, $\bar{m}\in {\rm l}^2(\mathbb{R}^2)$ are called \emph{equivalent from above} if the corresponding magnetic potentials $V$ and $\bar{V}$ (given by \eqref{eqn:mag2d}) are equal in the upper half-space, i.e., if $V\equiv \bar{V}$ in $\mathbb{R}^3_+$. They are called  \emph{equivalent from below} if $V$ and $\bar{V}$ are equal in the lower half-space, i.e., if $V\equiv \bar{V}$ in $\mathbb{R}^3_-$.  A magnetization $m$ is called \emph{silent from above} if it is equivalent from above to $\bar{m}\equiv0$ and \emph{silent from below} if it is equivalent from below to  $\bar{m}\equiv0$.
\end{defi}

A decomposition of $m$ that reflects this characterization is the so-called Hardy-Hodge decomposition (for details on the thin-plate case and all results mentioned in this section, the reader is referred to \cite{baratchart13}). For that purpose, we require the following vectorial operators:
\begin{align}
\bar{o}^{(1)}&=\left(R_1,R_2,\textnormal{id}\right),\label{eqn:eb1}
\\\bar{o}^{(2)}&=\left(-R_1,-R_2,\textnormal{id}\right),\label{eqn:eb2}
\\\bar{o}^{(3)}&=\left(-R_2,R_1,0\right),\label{eqn:eb3}
\end{align}
where $R_i$, $i=1,2$, are the Riesz transforms
\begin{align} 
R_i[F](x)=\lim_{\eps\to0}\frac{1}{2\pi}\int_{\mathbb{R}^2\setminus B_\varepsilon(x)}F(y)\frac{x_i-y_i}{|x-y|^3}{\dd}y,\qquad x=(x_1,x_2)\in\mathbb{R}^2,
\end{align}
of a scalar function $F$ of class $\LL^2(\mathbb{R}^2)$. We can now formulate the following theorem.

\begin{thm}[Hardy-Hodge Decomposition]\label{thm:hh2d}
Any function $f=(F_{x_1}, F_{x_2}, F_{x_3})\in {\rm l}^{2}(\mathbb{R}^2)$ can be decomposed into
\begin{align}\label{eqn:hh2ddecomp}
 f&=\bar{f}^{(1)}+\bar{f}^{(2)}+\bar{f}^{(3)}=\bar{o}^{(1)}[\bar{F}_1]+\bar{o}^{(2)}[\bar{F}_2]+\bar{o}^{(3)}[\bar{F}_3],
\end{align}
with scalar functions $\bar{F}_1$, $\bar{F}_2$, $\bar{F}_3$ given by
\begin{align}
 \bar{F}_1&=\frac{1}{2}\left(-R_1[F_{x_1}]-R_2[F_{x_2}]+F_{x_3}\right),\label{eqn:2df1}
 \\\bar{F}_2&=\frac{1}{2}\left(R_1[F_{x_1}]+R_2[F_{x_2}]+F_{x_3}\right),
 \\\bar{F}_3&=R_2[F_{x_1}]-R_1[F_{x_2}].\label{eqn:2df3}
\end{align}
\end{thm}

Spherical counterparts to the above theorem are introduced in more detail in Section \ref{sec:sdecomp}. A helpful formal notation for the operators $\bar{o}^{(1)}$, $\bar{o}^{(2)}$, $\bar{o}^{(3)}$ that emphasizes the connection between the spherical and the Euclidean thin-plate case is given in the next remark.

\begin{rem}\label{rem:not1}
It is well-known that the Riesz transform $R_i$ can formally also be expressed as $R_i=\nabla^*_i(\Delta^*)^{-\frac{1}{2}}$, $i=1,2$, where $\nabla^*_i$ denotes the $i$-th component of the surface gradient $\nabla^*$ from Section \ref{sec:notation}. With this notation at hand we can now reformulate the operators \eqref{eqn:eb1}--\eqref{eqn:eb3}:
\begin{align}
 \bar{o}^{(1)}&=o^{(1)}+o^{(2)}\left(\Delta^*\right)^{-\frac{1}{2}},\label{eqn:op2d1}
 \\\bar{o}^{(1)}&=o^{(1)}-o^{(2)}\left(\Delta^*\right)^{-\frac{1}{2}},\label{eqn:op2d2}
 \\\bar{o}^{(3)}&=o^{(3)}\left(\Delta^*\right)^{-\frac{1}{2}}.\label{eqn:op2d3}
\end{align}
\end{rem}

Supported by the Hardy-Hodge decomposition from Theorem \ref{thm:hh2d}, \cite{baratchart13} have derived several characterizations and uniqueness results under the constraint of unidirectionality and/or locally compact support on the magnetization $m$. We list those who relate to results for the spherical setting later on in Sections \ref{sec:sdecomp} and \ref{sec:sunique}.

\begin{thm}\label{thm:unque2d1}
 Let $m\in {\rm l}^2(\mathbb{R}^2)$ and $\bar{m}^{(1)}$, $\bar{m}^{(2)}$, $\bar{m}^{(3)}$ be given as in Theorem \ref{thm:hh2d}. Then the following assertions hold true:
 \begin{itemize}
  \item[(a)] The magnetization $\bar{m}^{(2)}$ is equivalent from above to $m$ while  $\bar{m}^{(1)}$ is equivalent from below to $m$.
  \item[(b)] The magnetization $m$ is silent from above if and only if $\bar{m}^{(2)}\equiv0$ while $m$ is silent from below if and only if $\bar{m}^{(1)}\equiv0$.
  \item[(c)] If $\textnormal{supp}(m)\subset \Gamma$, for a region $\Gamma\subset\mathbb{R}^2$ with $\Gamma\not=\mathbb{R}^2$, then $m$ is silent from above if and only if it is silent from below.
 \end{itemize}
\end{thm}


\begin{cor}[Unidirectional Magnetizations]\label{thm:unque2d3}
 Let $m\in {\rm l}^2(\mathbb{R}^2)$ be a non-tangential unidirectional magnetization, i.e., $m=Qv$ for a scalar function $Q\in {\rm L}^2(\mathbb{R}^2)$ and a fixed direction $v=(v_1,v_2,v_3)\in\mathbb{R}^3$, $v_3\not=0$. Furthermore, let $\Gamma\subset\mathbb{R}^2$ be a region with $\Gamma\not=\mathbb{R}^2$ and \emph{supp}$(m)\subset\Gamma$. 
 
Then $m$ is equivalent from above to no other unidirectional magnetization $\mathsf{m}=\mathsf{Q}\mathsf{v}$ with \emph{supp}$(\mathsf{m})\subset\Gamma$. Analogously,  $m$ is equivalent from below to no other unidirectional magnetization $\mathsf{m}=\mathsf{Q}\mathsf{v}$ with \emph{supp}$(\mathsf{m})\subset\Gamma$.
\end{cor}

\section{Spherical Decompositions}\label{sec:sdecomp}

From now on, we are strictly working in the spherical setting, i.e., we are investigating the magnetic potential $V$ that is generated by a magnetization $m\in {\rm l}^2(\Omega)$:
\begin{align}\label{eqn:smag0}
 V(x)=\frac{1}{4\pi}\int_{\Omega}m(\eta)\cdot\frac{x-\eta}{|x-\eta|^3}{\dd} \omega(\eta),\qquad x\in\mathbb{R}^3.
\end{align}

\subsection{Vector Spherical Harmonic Representation}\label{sec:vectharm}
A spherical version of the Hardy-Hodge decomposition from Theorem \ref{thm:hh2d} has been used in geomagnetic applications for quite some time in form of a decomposition of vector fields with respect to vector spherical harmonics $\tilde{y}^{(1)}_{n,k}$, $\tilde{y}^{(2)}_{n,k}$, $\tilde{y}^{(3)}_{n,k}$ (see, e.g., \cite{backus96, gubbins11, mayer06, mayermaier06, olsen10b}). These vector spherical harmonics can be defined via a suitable connection to the inner harmonics $H_{n,k}^{int}$ and the outer harmonics $H_{n,k}^{ext}$ (i.e., the harmonic extensions of scalar orthonormalized spherical harmonics $Y_{n,k}$ into $\Omega^{int}$ and $\Omega^{ext}$, respectively). More precisely,
\begin{align}
 \tilde{y}_{n,k}^{(1)}\left(\xi\right)&=-(\tilde{\mu}_n^{(1)})^{-\frac{1}{2}}\lim_{x\to\xi\atop x\in\Omega^{ext}}\nabla_xH_{n,k}^{ext}(x),&&n\in\mathbb{N}_0,\,k=1,\ldots 2n+1,\quad \xi\in\Omega,\label{eqn:ioh11}
 \\\tilde{y}_{n,k}^{(2)}\left(\xi\right)&=(\tilde{\mu}_n^{(2)})^{-\frac{1}{2}}\lim_{x\to\xi\atop x\in\Omega^{int}}\nabla_xH_{n,k}^{int}(x),&& n\in\mathbb{N},\,k=1,\ldots 2n+1, \quad \xi\in\Omega,\label{eqn:ioh22}
\end{align}
with normalization constants $\tilde{\mu}_{n}^{(1)}=(n+1)(2n+2)$, $\tilde{\mu}_{n}^{(2)}=n(2n+1)$. The vector spherical harmonics $\tilde{y}_{n,k}^{(3)}$, $n\in\mathbb{N},\,k=1,\ldots 2n+1$, are chosen such that, together with \eqref{eqn:ioh11} and \eqref{eqn:ioh22}, they form a complete orthonormal function systems in ${\rm l}^2(\Omega)$ with respect to the inner product $\langle\cdot,\cdot\rangle_{{\rm l}^2(\Omega)}$. These properties imply that a square-integrable vector field $f$ of the form $f=\nabla U$, where $U$ is harmonic in $\Omega^{ext}$, can be expressed by
\begin{align}
 f=\sum_{n=0}^\infty\sum_{k=1}^{2n+1} \langle f,\tilde{y}_{n,k}^{(1)}\rangle_{{\rm l}^2(\Omega)}\tilde{y}_{n,k}^{(1)}
\end{align}
on the sphere $\Omega$. Analogously,  a square-integrable  field $f$ of the form $f=\nabla U$, where $U$ is harmonic in $\Omega^{int}$, can be expressed by
\begin{align}
 f=\sum_{n=1}^\infty\sum_{k=1}^{2n+1} \langle f,\tilde{y}_{n,k}^{(2)}\rangle_{{\rm l}^2(\Omega)}\tilde{y}_{n,k}^{(2)}\label{eqn:fintrep}
\end{align}
on the sphere $\Omega$. More details on the involved (vector) spherical harmonics can be found, e.g., in \cite{freedenschreiner09}. Most general, one can state the following theorem.

\begin{thm}[Spherical Hardy-Hodge Decomposition I]\label{thm:shdh}
Any function $f\in {\rm l}^2(\Omega)$ can be decomposed into
\begin{align}
f&=\tilde{f}^{(1)}+\tilde{f}^{(2)}+\tilde{f}^{(3)}\label{eqn:split}
\\&=\sum_{n=0}^\infty\sum_{k=1}^{2n+1} \langle f,\tilde{y}_{n,k}^{(1)}\rangle_{{\rm l}^2(\Omega)}\tilde{y}_{n,k}^{(1)}+\sum_{n=1}^\infty\sum_{k=1}^{2n+1} \langle f,\tilde{y}_{n,k}^{(2)}\rangle_{{\rm l}^2(\Omega)}\tilde{y}_{n,k}^{(2)}+\sum_{n=1}^\infty\sum_{k=1}^{2n+1} \langle f,\tilde{y}_{n,k}^{(3)}\rangle_{{\rm l}^2(\Omega)}\tilde{y}_{n,k}^{(3)}.\nonumber
\end{align}
\end{thm}

\begin{rem}\label{rem:expdecomp}
 In order to investigate the consequences of Theorem \ref{thm:shdh} for the uniqueness of the magnetization $m$ and the magnetic potential $V$ in \eqref{eqn:smag0}, we first observe that, for $x\in\Omega^{ext}$ and $\eta\in\Omega$,
\begin{align}
\frac{x-\eta}{|x-\eta|^3}&=\lim_{y\to\eta\atop y\in\Omega^{int}}\nabla_y\frac{1}{|x-y|}.\label{eqn:1}
\end{align}
Clearly, $\frac{1}{|x-\cdot|}$ is harmonic in $\Omega^{int}$ and \eqref{eqn:fintrep} implies a representation of the form 
\begin{align}\label{eqn:fracrep}
\frac{x-\eta}{|x-\eta|^3}=\sum_{n=1}^\infty\sum_{k=1}^{2n+1}\left\langle \frac{x-\cdot}{|x-\cdot|^3},\tilde{y}_{n,k}^{(2)}\right\rangle_{{\rm l}^2(\Omega)}\tilde{y}_{n,k}^{(2)}(\eta), \quad x\in\Omega^{ext},\eta\in\Omega.
\end{align}
In detail, using the addition theorem for scalar spherical harmonics and a series representation of $\frac{1}{|x-\cdot|}$ in terms of Legendre polynomials $P_n$, we obtain
\begin{align}
\frac{x-\eta}{|x-\eta|^3}&=\lim_{y\to\eta\atop y\in\Omega^{int}}\nabla_y\frac{1}{|x|}\sum_{n=0}^\infty\left(\frac{|y|}{|x|}\right)^nP_n\left(\frac{x}{|x|}\cdot\frac{y}{|y|}\right)\label{eqn:rep}
\\&=\lim_{y\to\eta\atop y\in\Omega^{int}}\sum_{n=1}^\infty\frac{4\pi}{2n+1}\sum_{k=1}^{2n+1}H_{n,k}^{ext}\left(x\right) \nabla_yH_{n,k}^{int}\left(y\right)\nonumber
\\&=\sum_{n=1}^\infty\sum_{k=1}^{2n+1}\left(\frac{4\pi(\tilde{\mu}_{n}^{(2)})^{\frac{1}{2}}}{2n+1}H_{n,k}^{ext}\left(x\right) \right)\tilde{y}_{n,k}^{(2)}\left(\eta\right).\nonumber
\end{align}
If we now substitute \eqref{eqn:fracrep} or \eqref{eqn:rep} into the representation \eqref{eqn:smag0} of $V$, we see that, due to the orthonormality of the vector spherical harmonics, the contribution $\tilde{m}^{(2)}$ of $m=\tilde{m}^{(1)}+\tilde{m}^{(2)}+\tilde{m}^{(3)}$ generates the exact same magnetic potential $V$ in $\Omega^{ext}$ as the entire magnetization $m$ itself. Computations analogous to \eqref{eqn:1}, \eqref{eqn:fracrep}, and \eqref{eqn:rep} for the case $x\in\Omega^{int}$ and $\eta\in\Omega$ lead to the conclusion that the contribution $\tilde{m}^{(1)}$ of $m=\tilde{m}^{(1)}+\tilde{m}^{(2)}+\tilde{m}^{(3)}$ generates the exact same magnetic potential $V$ in $\Omega^{int}$ as the entire magnetization $m$ itself. 
\end{rem}

From Remark \ref{rem:expdecomp}, it becomes clear that the spherical decomposition of Theorem \ref{thm:shdh} reveals the desired properties corresponding to the thin-plate case in Theorem \ref{thm:unque2d1}(a),(b) (which has already been observed in \cite{gubbins11}).

\subsection{Operator Representation}\label{sec:op}
In order to be able to obtain a spherical version of Theorem \ref{thm:unque2d1}(c), we reformulate the decomposition of Theorem \ref{thm:shdh} in terms of a set of pseudo-differential operators $\oox$, $\ooy$, $\ooz$  as  indicated in \cite{freedengerhards12, freedenschreiner09, gerhards11a, gerhards12}. More precisely,
\begin{align}
 \oox&=\ox\left(\DD+\frac{1}{2}\right)-\oy,\label{eqn:oo1}
 \\\ooy&=\ox\left(\DD-\frac{1}{2}\right)+\oy,\label{eqn:oo2}
 \\\ooz&=\oz,\label{eqn:oo3}
\end{align}
where $\DD$ denotes the pseudo-differential operator
\begin{align}\label{eqn:d}
  \DD=\left(-\Delta^*+\frac{1}{4}\right)^{\frac{1}{2}}.
\end{align}
The previously introduced vector spherical harmonics can then be alternatively expressed by
\begin{align}
 \tilde{y}_{n,k}^{(1)}&=(\tilde{\mu}_{n}^{(1)})^{-\frac{1}{2}}\oox Y_{n,k},\qquad n\in\mathbb{N}_0,\,k=1,\ldots 2n+1,\label{eqn:vh1}
 \\\tilde{y}_{n,k}^{(2)}&=(\tilde{\mu}_{n}^{(2)})^{-\frac{1}{2}}\ooy Y_{n,k},\qquad n\in\mathbb{N},\,k=1,\ldots 2n+1,\label{eqn:vh2}
 \\\tilde{y}_{n,k}^{(3)}&=(\tilde{\mu}_{n}^{(3)})^{-\frac{1}{2}}\ooz Y_{n,k},\qquad n\in\mathbb{N},\,k=1,\ldots 2n+1,\label{eqn:vh3}
\end{align}
with normalization constants $\tilde{\mu}_{n}^{(1)}=(n+1)(2n+2)$, $\tilde{\mu}_{n}^{(2)}=n(2n+1)$, and $\tilde{\mu}_{n}^{(3)}=n(n+1)$. In other words, the properties of the decomposition from Theorem \ref{thm:shdh} carry over to a decomposition with respect to the operators  $\oox$, $\ooy$, $\ooz$. To formulate such a decomposition, we first recapitulate the spherical Helmholtz decomposition (see, e.g., \cite{freedengerhards12, freedenschreiner09} and references therein). 

\begin{thm}[Spherical Helmholtz Decomposition]\label{thm:helmdecomp}
Any function $f\in {\rm l}^2(\Omega)$ can be decomposed into
 \begin{align}
  f=\ox [F_1]+\oy [F_2]+\oz [F_3],
 \end{align}
where the scalar functions $F_1$, $F_2$, $F_3$ are uniquely determined by the conditions $\int_\Omega F_2(\eta) {\dd}\omega(\eta)=\int_\Omega F_3(\eta) {\dd}\omega(\eta)=0$. 
\end{thm}

The spherical Helmholtz decomposition simply states a decomposition of a spherical vector field into a normal component and two tangential components (of which one is surface divergence-free and the other one surface curl-free). From \cite{freedengerhards12, gerhards11a,gerhards12}, we now take the following theorem.

\begin{thm}[Spherical Hardy-Hodge Decomposition II]\label{thm:shhdecomp}
Any function $f\in {\rm l}^2(\Omega)$ can be decomposed into
 \begin{align}\label{eqn:shhdecomp}
  f=\fx+\fy+\fz=\oox [\Fx]+\ooy [\Fy]+\ooz [\Fz],
 \end{align}
where the scalar functions $\Fx$, $\Fy$, $\Fz$ are uniquely determined by the conditions $\int_\Omega \Fx(\eta)-\Fy(\eta) {\dd}\omega(\eta)=\int_\Omega \Fz(\eta) {\dd}\omega(\eta)=0$. If $F_1$, $F_2$, $F_3$ are the Helmholtz scalars of $f$ as given in Theorem \ref{thm:helmdecomp}, then
\begin{align}
 \Fx&=\frac{1}{2}\left(\DD^{-1}[F_1]-F_2+\frac{1}{2}\DD^{-1}[F_2]\right),\label{eqn:fx}
 \\\Fy&=\frac{1}{2}\left(\DD^{-1}[F_1]+F_2+\frac{1}{2}\DD^{-1}[F_2]\right),\label{eqn:fy}
 \\\Fz&=F_3.\label{eqn:fz}
\end{align}
\end{thm}

A slight modification of the operators $\oox$, $\ooy$, $\ooz$ highlights the relation of Theorem \ref{thm:shhdecomp} to the Euclidean thin-plate case in Theorem \ref{thm:hh2d} and Remark \ref{rem:not1}. Application of the operator $\left(\DD+\frac{1}{2}\right)^{-1}$ to $\oox$ and application of $\left(\DD-\frac{1}{2}\right)^{-1}$ to $\ooy$ and $\ooz$ leads to the operators 
\begin{align}
 \ooox&=\ox-\oy\left(\DD+\frac{1}{2}\right)^{-1},\label{eqn:ooo1}
 \\\oooy&=\ox+\oy\left(\DD-\frac{1}{2}\right)^{-1},\label{eqn:ooo2}
 \\\oooz&=\oz\left(\DD-\frac{1}{2}\right)^{-1}.\label{eqn:ooo3}
\end{align}
Comparing this to \eqref{eqn:eb1}--\eqref{eqn:eb3} and \eqref{eqn:op2d1}--\eqref{eqn:op2d3}, we see that $o^{(2)}\left(\DD+\frac{1}{2}\right)^{-1}$ and $o^{(2)}\left(\DD-\frac{1}{2}\right)^{-1}$, respectively, take the role of the Riesz transform in the thin-plate case. However, it should be remarked that in the literature the Riesz transform on the sphere is typically given by $\oy(\Delta^*)^{-\frac{1}{2}}$ (see, e.g., \cite{daixu13}). Furthermore, $\left(\DD-\frac{1}{2}\right)^{-1}$ is well-defined only if restricted to the space $\LL_0^2(\Omega)=\{F\in \LL^2(\Omega):\int_\Omega F(\eta){\dd}\omega(\eta)=0\}$ since the constant functions form the nullspace of $\DD-\frac{1}{2}$. Yet, this is no restriction for our further considerations. As an alternative to Theorem \ref{thm:shhdecomp}, we can now state the following theorem.
 
\begin{thm}[Spherical Hardy-Hodge Decomposition III]\label{thm:shhdecomp2}
Any function $f\in {\rm l}^2(\Omega)$ can be decomposed into
 \begin{align}\label{eqn:shhdecomp2}
  f=\bar{f}_1+\bar{f}_2+\bar{f}_3=\ooox [\Ffx]+\oooy [\Ffy]+\oooz [\Ffz],
 \end{align}
where the scalar functions $\Ffx$, $\Ffy$, $\Ffz$ are uniquely determined by the conditions $\int_\Omega \Ffy(\eta) {\dd}\omega(\eta)=\int_\Omega \Ffz(\eta) {\dd}\omega(\eta)=0$. If $F_1$, $F_2$, $F_3$ are the Helmholtz scalars of $f$ as given in Theorem \ref{thm:helmdecomp}, then
\begin{align}
 \Ffx&=\frac{1}{2}\left(F_1+\frac{1}{2}\DD^{-1}[F_1]-\DD[F_2]+\frac{1}{4}\DD^{-1}[F_2]\right),\label{eqn:ffx}
 \\\Ffy&=\frac{1}{2}\left(F_1-\frac{1}{2}\DD^{-1}[F_1]+\DD[F_2]-\frac{1}{4}\DD^{-1}[F_2]\right),\label{eqn:ffy}
 \\\Ffz&=\DD[F_3]-\frac{1}{2}F_3.\label{eqn:ffz}
\end{align}
\end{thm}

\begin{prf}
For $\Ffx=(\DD+\frac{1}{2})[\Fx]$, $\Ffy=(\DD-\frac{1}{2})[\Fy]$, and  $\Ffz=(\DD-\frac{1}{2})[\Fz]$, we directly obtain the representations \eqref{eqn:ffx}--\eqref{eqn:ffz} from the corresponding representations of $\Fx$, $\Fy$, and $\Fz$ in Theorem \ref{thm:shhdecomp}. Concerning the uniqueness, we can restrict the considerations to the case $f\equiv 0$. Using \eqref{eqn:ooo1}--\eqref{eqn:ooo2} in \eqref{eqn:shhdecomp2} leads to 
\begin{align}
0&\equiv\ooox [\Ffx]+\oooy [\Ffy]+\oooz [\Ffz]
\\&=o^{(1)}\left(\Ffx+\Ffy\right)+o^{(2)}\left(- \left(\DD+\frac{1}{2}\right)^{-1}[\Ffx]+\left(\DD-\frac{1}{2}\right)^{-1}[\Ffy]\right)+o^{(3)}\Ffz.
\end{align}
The uniqueness of the normal part of $f$ leads to $\Ffx=-\Ffy$. Thus, the assumption $\int_\Omega \Ffy(\eta) {\dd}\omega(\eta)=0$ implies $\int_\Omega \Ffx(\eta) {\dd}\omega(\eta)=0$ and, consequently,
\begin{align}
\int_\Omega -\left(\DD+\frac{1}{2}\right)^{-1}[\Ffx](\eta)+\left(\DD-\frac{1}{2}\right)^{-1}[\Ffy](\eta)\,{\dd}\omega(\eta)=0.
\end{align}
The uniqueness of the spherical Helmholtz decomposition in Theorem \ref{thm:helmdecomp} then leads to $(\DD+\frac{1}{2})^{-1}[\Ffx]+(\DD-\frac{1}{2})^{-1}[\Ffx]=(\DD+\frac{1}{2})^{-1}[\Ffx]-(\DD-\frac{1}{2})^{-1}[\Ffy]\equiv0$, or in other words, by application of $\DD-\frac{1}{2}$,
\begin {align}
\left(2-\left(\DD+\frac{1}{2}\right)^{-1}\right)[\Ffx]\equiv0.
\end {align}
Since $2-\left(\DD+\frac{1}{2}\right)^{-1}$ is injective on $\LL_0^2(\Omega)$, we get $\Ffx=-\Ffy\equiv0$. Moreover, the condition $\int_\Omega \Ffz(\eta) {\dd}\omega(\eta)=0$ and the uniqueness of the Helmholtz scalars in Theorem \ref{thm:helmdecomp} implies $\Ffz\equiv0$ and concludes the proof.
\end{prf}

\begin{rem}
 We need to emphasize that the functions $\fx$, $\fy$, $\fz$ in Theorems \ref{thm:shdh} and \ref{thm:shhdecomp} and the functions $\ffx$, $\ffy$, $\ffz$ in Theorem \ref{thm:shhdecomp2} are identical. The theorems only differ in the applied operators and the corresponding scalar functions $\Fx$, $\Fy$, $\Fz$, and $\Ffx$, $\Ffy$, $\Ffz$. For our later considerations, we work with the operators $\oox$, $\ooy$, $\ooz$ and the functions $\Fx$, $\Fy$, $\Fz$.
\end{rem}

\section{Uniqueness Issues for Spherical Magnetizations}\label{sec:sunique}

We start by defining the notion of equivalence of magnetizations, analogous to the thin-plate case  in Definition \ref{def:equivabove}.

\begin{defi}
Two magnetizations $m$, $\bar{m}\in {\rm l}^2(\Omega)$ are called \emph{equivalent from inside} if the corresponding magnetic potentials $V$ and $\bar{V}$ (given by \eqref{eqn:smag0}) are equal in the interior $\Omega^{int}$, i.e., if $V\equiv \bar{V}$ in $\Omega^{int}$. They are called  \emph{equivalent from outside} if $V$ and $\bar{V}$ are equal in the exterior $\Omega^{ext}$, i.e., if $V\equiv \bar{V}$ in $\Omega^{ext}$.  A magnetization $m$ is called \emph{silent from inside} if it is equivalent from inside to $\bar{m}\equiv0$ and \emph{silent from outside} if it is equivalent from outside to  $\bar{m}\equiv0$.
\end{defi}
Theorem \ref{thm:shhdecomp} eventually allows the following characterizations of spherical magnetizations. 

\begin{thm}\label{thm:sunique1}
 Let $m\in {\rm l}^2(\Omega)$ and $\tilde{m}^{(1)}$, $\tilde{m}^{(2)}$, $\tilde{m}^{(3)}$ be given as in Theorem \ref{thm:shhdecomp}. Then the following assertions hold true:
 \begin{itemize}
  \item[(a)] The magnetization $\tilde{m}^{(2)}$ is equivalent from outside to $m$ while  $\tilde{m}^{(1)}$ is equivalent from inside to $m$.
  \item[(b)] The magnetization $m$ is silent from outside if and only if $\tilde{m}^{(2)}\equiv0$ while $m$ is silent from inside if and only if $\tilde{m}^{(1)}\equiv0$.
  \item[(c)] If $\textnormal{supp}(m)\subset \Gamma$, for a region $\Gamma\subset\Omega$ with $\Gamma\not=\Omega$, then $m$ is silent from outside if and only if it is silent from inside.
  \end{itemize}
\end{thm} 

\begin{prf}
Parts (a) and (b) are direct consequences of the considerations in Remark \ref{rem:expdecomp}. Concerning part (c), we assume that $m$ is silent from outside and can now conclude that $\tilde{m}^{(2)}\equiv0$ and, therefore,
\begin{align}\label{eqn:tm1zero}
 \tilde{M}_2=\frac{1}{2}\left(\DD^{-1}[M_1]+M_2+\frac{1}{2}\DD^{-1}[M_2]\right)=0,
\end{align}
where $M_1$, $M_2$, $M_3$ are the Helmholtz scalars of $m$ as supplied in Theorem \ref{thm:helmdecomp}. Plugging \eqref{eqn:tm1zero} into \eqref{eqn:fx} implies $\tilde{M}_1=-M_2$. In other words,
\begin{align}\label{eqn:msplit}
 m=\tilde{m}^{(1)}+\tilde{m}^{(3)}=-\oox M_2+\ooz M_3.
\end{align}
Expanding the Helmholtz scalar $M_2$ in terms of spherical harmonics leads to the expression $\oox M_2=\sum_{n=0}^\infty\sum_{k=1}^{2n+1}\langle M_2,Y_{n,k}\rangle_{L^2(\Omega)}\oox Y_{n,k}$. Next, for $x\in\Omega^{ext}$, we set
\begin{align}
 U(x)&=\sum_{n=0}^\infty\sum_{k=1}^{2n+1}\langle M_2,Y_{n,k}\rangle_{L^2(\Omega)}H_{n,k}^{ext}\left(x\right)\label{eqn:u}
 \\&=\int_\Omega \left(\sum_{n=0}^\infty\sum_{k=1}^{2n+1}H_{n,k}^{ext}\left(x\right)Y_{n,k}\left(\eta\right)\right)M_2(\eta){\dd}\omega(\eta)\nonumber
 \\&=\int_\Omega \left(\sum_{n=0}^\infty\frac{2n+1}{4\pi}|x|^{-(n+1)}P_n\left(\frac{x}{|x|}\cdot\eta\right)\right)M_2(\eta){\dd}\omega(\eta)\nonumber
 \\&=-\frac{1}{4\pi}\int_\Omega \frac{1-|x|^2}{\left(1+|x|^2-2x\cdot \eta\right)^{\frac{3}{2}}}M_2(\eta){\dd}\omega(\eta).\nonumber
\end{align}
The last equalities can be found, e.g., in \cite{freedenschreiner09} and references therein. Clearly, $U$ is harmonic in $\Omega^{ext}$. Remembering \eqref{eqn:ioh11}, \eqref{eqn:ioh22}, and \eqref{eqn:vh1}--\eqref{eqn:vh3}, we can conclude that 
\begin{align}\label{eqn:nablaurep}
 \lim_{x\to\xi\atop x\in\Omega^{ext}}\nabla_x U(x)=-\oox_\xi M_2(\xi),\qquad \xi\in\Omega.
\end{align}
From \eqref{eqn:tm1zero} we find $M_2=-(\DD+\frac{1}{2})^{-1}M_1$, and because $(\DD+\frac{1}{2})^{-1}$ is selfadjoint, it holds 
\begin{align}
 U(x)&=\frac{1}{4\pi}\int_\Omega \left(\left(\DD_\eta+\frac{1}{2}\right)^{-1}\frac{1-|x|^2}{\left(1+|x|^2-2x\cdot \eta\right)^{\frac{3}{2}}}\right)M_1(\eta){\dd}\omega(\eta).
\end{align}
Since $\textnormal{supp}(m)\subset\Gamma$, we get $M_1\equiv0$ in $\Gamma^c=\Omega\setminus\overline{\Gamma}$ and can extend $U$ across $\Gamma^c$ to a harmonic function in the cone $\mathcal{C}_{\Gamma^c}=\{x\in\mathbb{R}^3\setminus\{0\}:\frac{x}{|x|}\in\Gamma^c\}$.

Next, we observe that there must exist some function $\bar{W}$ on ${\Gamma^c}$ such that $\tilde{m}^{(3)}=o^{(2)}\bar{W}$ in ${\Gamma^c}$. This follows from \eqref{eqn:msplit} because $m\equiv0$ in ${\Gamma^c}$, because $\tilde{m}^{(3)}$ is tangential, and because $\tilde{m}^{(1)}$ is solely composed by a normal part due to $o^{(1)}$ and a tangential part due to $o^{(2)}$. Now, setting $W(x)=\bar{W}(x/|x|)$, for $x\in\mathcal{C}_{\Gamma^c}$, we see that $W$ is harmonic in $\mathcal{C}_{\Gamma^c}$. Combining our findings up to now, we have a function $U+W$ that is harmonic in $\mathcal{C}_{\Gamma^c}$ and which satisfies
\begin{align}
 m(\xi)=\lim_{x\to\xi\atop x\in\Omega^{ext}}\nabla_x(U(x)+W(x))=0,\qquad \xi\in{\Gamma^c}.
\end{align}
Therefore, the normal and the tangential derivative of $U+W$ vanish on ${\Gamma^c}$. Consequently, since $U+W$ is harmonic in $\mathcal{C}_{\Gamma^c}$, typical analytical continuation arguments lead to $U+W\equiv0$ in $\mathcal{C}_{\Gamma^c}$. In particular, $\frac{x}{|x|}\cdot\nabla_x(U(x)+W(x))=\frac{x}{|x|}\cdot\nabla_xU(x)=0$, for $x\in\mathcal{C}_{\Gamma^c}$. Since $\nabla U$ is harmonic in $\Omega^{ext}$ by the construction in \eqref{eqn:u}, we obtain $\frac{x}{|x|}\cdot\nabla_xU(x)=0$, for $x\in\Omega^{ext}$. Combining this with \eqref{eqn:nablaurep} leads to 
\begin{align}
\lim_{x\to\xi\atop x\in\Omega^{ext}}\frac{x}{|x|}\cdot \nabla_x U(x)=-\xi\cdot\oox_\xi M_2(\xi)=-\left(D+\frac{1}{2}\right)[M_2](\xi)=0,\qquad \xi\in\Omega.
\end{align}
In other words $M_2\equiv0$. Observing \eqref{eqn:tm1zero}, we now find $\DD^{-1}[M_1]\equiv0$, which leads to $M_1\equiv0$ and eventually to $\tilde{m}^{(1)}\equiv0$ by \eqref{eqn:fx}. Thus,
\begin{align}
 m=\tilde{m}^{(3)},
\end{align}
which, by parts (a) and (b), implies that $m$ is silent from outside as well as inside. 
\end{prf}

Theorem \ref{thm:sunique1} has some direct consequences for the uniqueness of induced magnetizations of the form 
\begin{align}\label{eqn:indmag}
 m(\xi)=Q(\xi)v(\xi),\quad\xi\in\Omega,
\end{align}
where $v$ is a known inducing vector field and $Q$ an unknown scalar susceptibility. If supp($m$)$\subset\Gamma$ and if there exists another magnetization $\mathsf{m}=\mathsf{Q}v$ with supp($\mathsf{m}$)$\subset\Gamma$ that generates the same magnetic potential in $\Omega^{ext}$, then $m-\mathsf{m}$ is silent from outside and Theorem \ref{thm:sunique1} implies that $m-\mathsf{m}\equiv \ooz {M}$ for some scalar function ${M}$. In particular, $m-\mathsf{m}=(Q-\mathsf{Q})v$ must be tangential. If the known inducing field $v$ is non-tangential (i.e., $v(\xi)\cdot\xi\not=0$ for some $\xi\in\Omega$), this implies $m=\mathsf{m}$. Put in more rigorous terms, we obtain the following corollary. 

\begin{cor}[Uniqueness of Induced Magnetizations]\label{thm:indmag}
Let $m\in {\rm l}^2(\Omega)$ be of the form $m=Qv$ (where $Q\in\LL^2(\Omega)$ is unknown and $v\in{\rm l}^\infty(\Omega)$ is known) with \emph{supp}$(m)\subset\Gamma$ for a fixed region $\Gamma\subset\Omega$ with $\Gamma\not=\Omega$. Furthermore, we assume that $v$ is non-tangential.
 
Then $m$ is equivalent from outside to no other induced magnetization $\mathsf{m}=\mathsf{Q}v\in {\rm l}^2(\Omega)$  with \emph{supp}$(\mathsf{m})\subset\Gamma$. Analogously, $m$ is equivalent from inside to no other induced magnetization  $\mathsf{m}=\mathsf{Q}v\in {\rm l}^2(\Omega)$ with \emph{supp}$(\mathsf{m})\subset\Gamma$.
\end{cor}

The corollary above suits an application to the induced Earth's crustal magnetization since the inducing field $v$ is typically the Earth's core magnetic field, which is non-tangential and known from models such as \cite{olsen14}. A restriction to unidirectional magnetizations as in Corollary \ref{thm:unque2d3} (which demands further structure on the inducing field $v$ but, in exchange, does not require to know $v$ in advance)  for the Euclidean thin-plate case has no particular relevance in a global spherical context. Actually, depending on how one defines spherical unidrectionality, a spherical counterpart to Corollary \ref{thm:unque2d3} does not necessarily hold true. 
 
 \begin{figure}
\begin{center}
\includegraphics[scale=0.21]{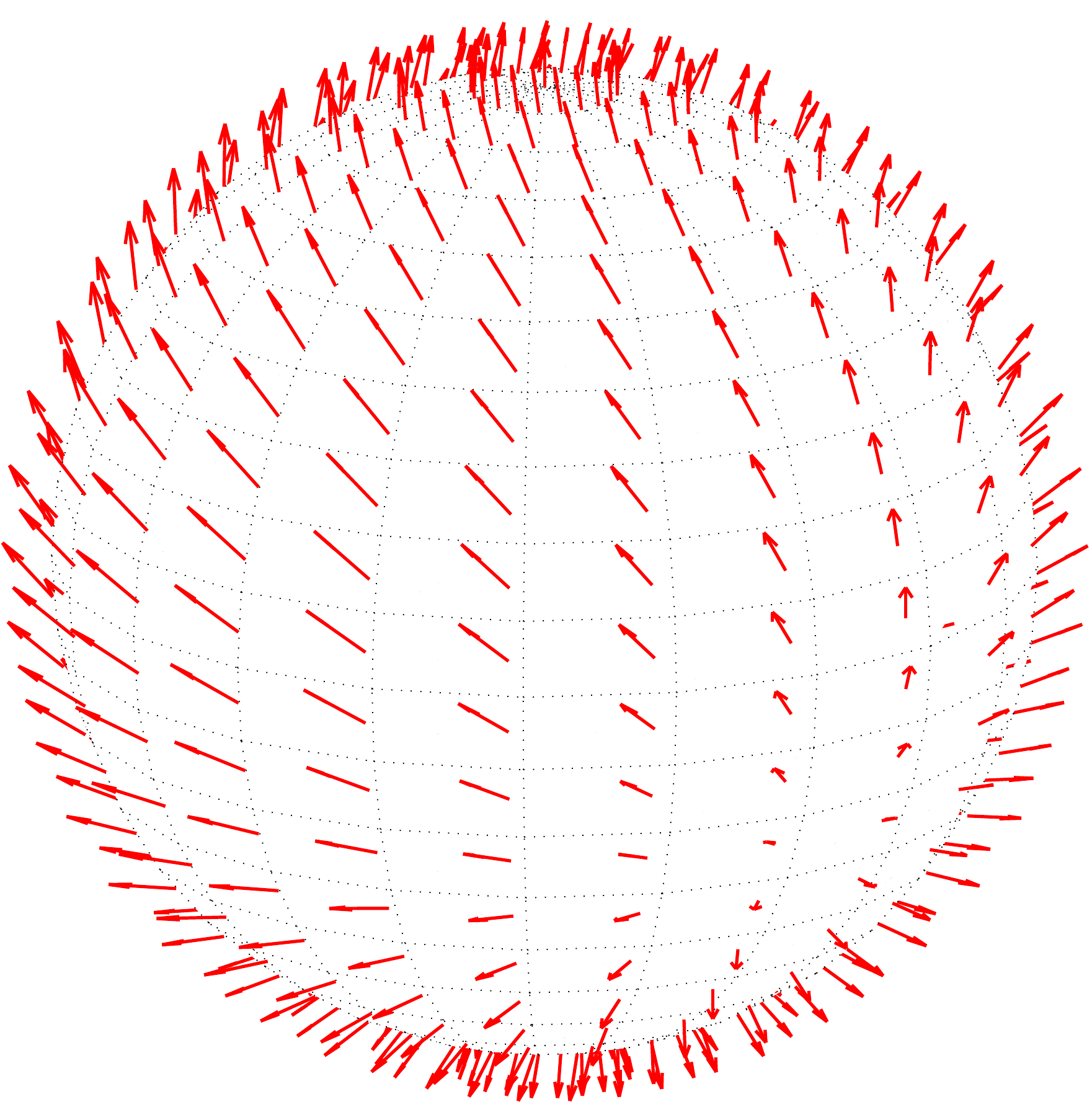}\qquad\qquad\qquad\includegraphics[scale=0.21]{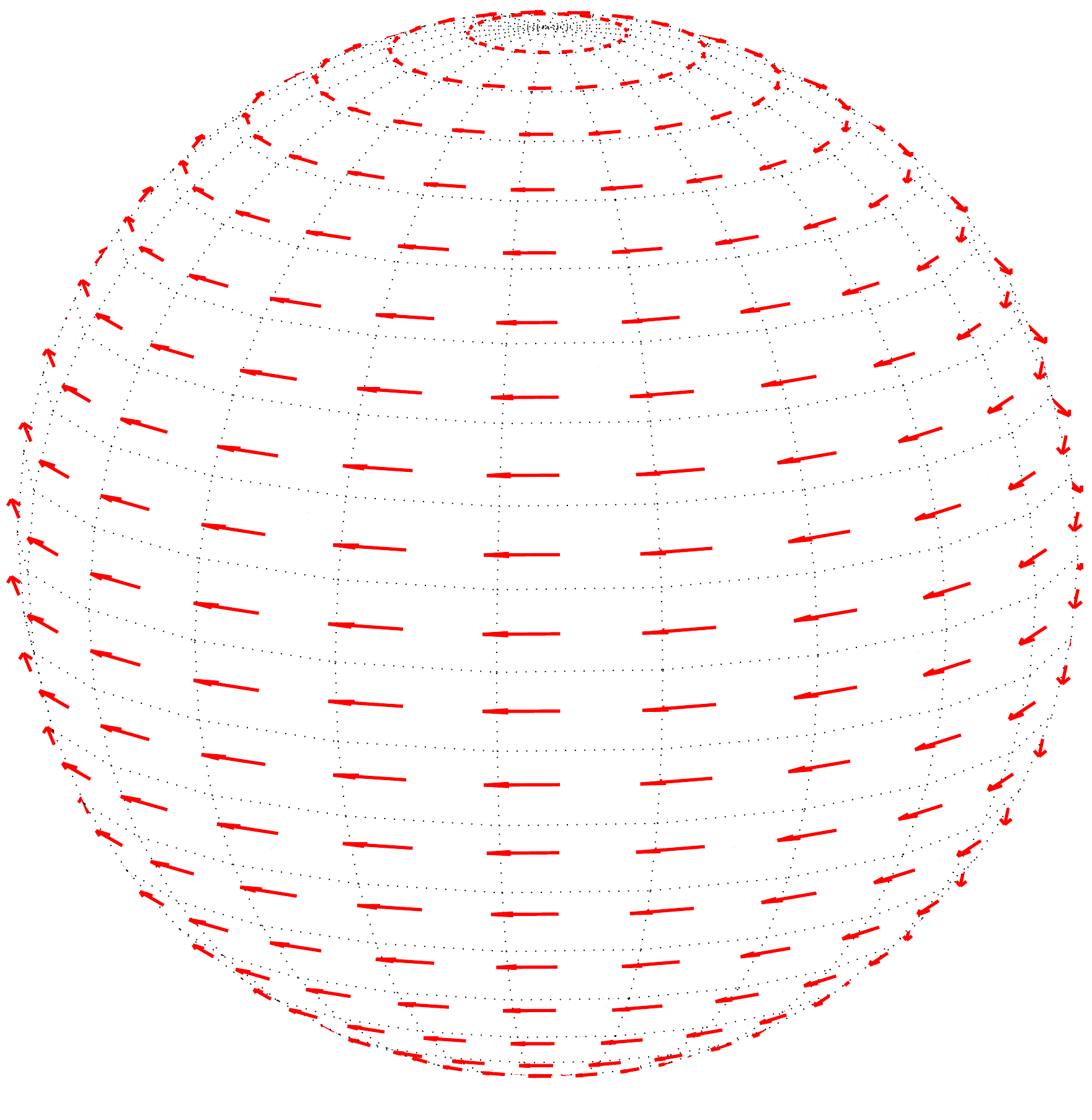}
\end{center}
\caption{Example of vector fields $m$ according to \eqref{eqn:sunid}, with $Q\equiv1$, $\zeta=(0,0,1)$, and $v_1=v_2=v_3=1$ (left) and  with $Q\equiv1$, $\zeta=(0,0,1)$, and $v_1=v_2=0$, $v_3=1$ (right).}
\label{fig:unidic}
\end{figure}

\begin{rem}
A fairly natural definition is to call a spherical magnetization unidirectional if it is constant with respect to a certain local spherical coordinate system. For simplicity, we regard magnetizations of the form
\begin{align}\label{eqn:sunid}
 m(\xi)=Q(\xi)\left(v_1\xi+v_2 (\zeta-(\zeta\cdot\xi)\xi)+v_3\xi\wedge\zeta\right), \quad\xi\in\Omega,
\end{align}
where $\zeta\in\Omega$ and $v_1,v_2,v_3\in\mathbb{R}$ are fixed and $Q\in \LL^2(\Omega)$ (an illustration is given in Figure \ref{fig:unidic}). Note that the tangential components $\zeta-(\zeta\cdot\xi)\xi$ and $\xi\wedge\zeta$ vanish as $\xi$ approaches the poles $\pm\zeta$. Therefore, the vectors $\xi$, $\zeta-(\zeta\cdot\xi)\xi$, and $\xi\wedge\zeta$ do not form a local coordinate system in the points $\xi=\pm\zeta$. Yet, the example derived below can serve as a counterexample to uniqueness for spherical unidirectional magnetizations if we assume that the function $Q$ vanishes around $\pm\zeta$.

A consequence of a spherical counterpart to Corollary \ref{thm:unque2d3} would be that the trivial magnetization $m\equiv0$ is the only spherical unidirectional magnetization that is silent from outside and has support supp$(m)\subset\Gamma$. This particularly includes purely tangential magnetizations
\begin{align}
 m(\xi)=Q(\xi)\left(v_2 (\zeta-(\zeta\cdot\xi)\xi)+v_3\xi\wedge\zeta\right), \quad\xi\in\Omega.
\end{align}
If such a magnetization is silent from outside and supp$(m)\subset\Gamma$, then the proof of Theorem \ref{thm:sunique1} implies that $M_2\equiv0$, where $M_2$ is the Helmholtz scalar of $m$ as indicated in Theorem \ref{thm:helmdecomp}. Applying the Beltrami operator and observing the properties of $M_2$ leads to
\begin{align}
 0&=\Delta^*M_2(\xi)=\nabla^*\cdot m(\xi)\label{eqn:dm20}
 \\&=\nabla^*Q(\xi)\cdot\left(v_2(\zeta-(\xi\cdot\zeta)\xi)+v_3\xi\wedge\zeta\right)+Q(\xi)\left(-2v_2(\xi\cdot\zeta)+v_2\nabla^*_\xi\cdot(\xi\wedge\zeta)\right)\nonumber
 \\&=\zeta\cdot\left(v_2\nabla^*Q(\xi)-v_3\LL^*Q(\xi)\right)-2v_2(\zeta\cdot \xi)Q(\xi).\nonumber
\end{align}
We now choose $\zeta=(0,0,1)$. Parametrizing the sphere $\Omega$ with respect to polar distance $t\in[-1,1]$ and longitude $\varphi\in[0,2\pi)$, and using representations of $\nabla^*$ and $\LL^*$ with respect to this parametrization (see, e.g., \cite{freedengerhards12,freedenschreiner09}), equation \eqref{eqn:dm20} can be rewritten as follows:
\begin{align}
 0=v_2(1-t^2)\frac{\partial}{\partial t}Q(t,\varphi)-v_3 \frac{\partial}{\partial \varphi}Q(t,\varphi)-2v_2Q(t,\varphi).
\end{align}
We further simplify the above by choosing $v_2=0$ and obtain
\begin{align}
 0=v_3 \frac{\partial}{\partial \varphi}Q(t,\varphi),
\end{align}
implying $Q(t,\varphi)=\mathsf{Q}(t)$, for some function $\mathsf{Q}$ that only depends on the polar distance $t$. Summarizing the above considerations, we can conclude that a function $\mathsf{Q}\in \CC^{(2)}([-1,1])$ with  $\mathsf{Q}(t)=0$, for $t=\xi\cdot\zeta$ and $\xi\in\Gamma^c=\Omega\setminus\overline{\Gamma}$, yields a non-trivial spherical unidirectional magnetization
\begin{align}
 m(\xi)=v_3\mathsf{Q}(\xi\cdot\zeta)\xi\wedge\zeta, \quad \xi\in\Omega,
\end{align}
that is silent from outside and satisfies supp$(m)\subset\Gamma$. In other words, we cannot expect to obtain a spherical counterpart to Corollary \ref{thm:unque2d3} for spherical unidirectionality as meant above.

On the other hand, if we call a spherical magnetization unidirectional if $m(\xi)=Q(\xi)v$, for a fixed $v=(v_1,v_2,v_3)\in\mathbb{R}^3\setminus\{0\}$ (i.e., if we use an identical definition of unidirectionality as in the planar case of Corollary \ref{thm:unque2d3}), we were not able to confirm nor to disprove uniqueness. 
\end{rem}

\begin{rem}\label{rem:locexist}
The existence of an induced magnetization $\mathsf{m}=Qv$ with supp$(\mathsf{m})\subset\Gamma$ that is equivalent from outside to a (not necessarily induced) magnetization $m\in{\rm l}^2(\Omega)$ with supp$(m)\subset\Gamma$ is generally not guaranteed. Assume, e.g., that $v(\xi)=\xi$, $\xi\in\Omega$.  If $m$ and $\mathsf{m}$ are equivalent from outside, then $m-\mathsf{m}$ is silent from outside and supp$(m-\mathsf{m})\subset\Gamma$. The proof of Theorem \ref{thm:sunique1} yields that this holds true if and only if $M_1=\mathsf{M}_1$ and $M_2=\mathsf{M}_2$, where $M_1$, $M_2$, and $\mathsf{M}_1$, $\mathsf{M}_2$ are the Helmholtz scalars of $m$ and $\mathsf{m}$, respectively. However, for our choice of $v$ it is clear that $\mathsf{M}_2\equiv0$. In other words, for any magnetization $m$ with $M_2\not\equiv0$ and supp$(m)\subset\Gamma$ there does not exist an induced magnetization $\mathsf{m}=Qv$ with supp$(\mathsf{m})\subset\Gamma$ that is equivalent from outside to $m$.
 \end{rem}

The situation of existence of an induced magnetization changes if we drop the condition that $\mathsf{m}=Qv$ has to satisfy supp$(\mathsf{m})\subset\Gamma$, at least if the inducing field $v$ satisfies certain conditions.

\begin{defi}\label{def:adm}
 We call a vector field $v\in{\rm l}^\infty(\Omega)$ \emph{admissible} if $|\xi\cdot v(\xi)|\geq C$, $\xi\in\Omega$, for some constant $C>0$, and if the coefficients
 \begin{align}
  v_{n,k,m,l}=\int_\Omega Y_{m,l}(\eta)\frac{v(\eta)}{\eta\cdot v(\eta)}\cdot \nabla^*Y_{n,k}(\eta)d\omega(\eta)
  \end{align}
  satisfy the following property: For every $n^3$-summable sequence $g_{n,k}$, $n\in\mathbb{N}_0$, $k=1,\ldots,2n+1$, i.e., $\sum_{n=0}^\infty\sum_{k=1}^{2n+1}n^6|g_{n,k}|^2<\infty$, the infinite dimensional system of linear equations
\begin{align}
 \sum_{m=0}^\infty\sum_{l=1}^{2m+1}\left(\frac{1}{nn}v_{n,k,m,l}+\frac{1}{n+\frac{1}{2}}\delta_{n,m}\delta_{k,l}\right)\gamma_{m,l}=&g_{n,k},\quad n\in\mathbb{N}_0,k=1,\ldots,2n+1,\label{eqn:inflin}
\end{align}
has a $n^2$-summable solution $\gamma_{n,k}$, $n\in\mathbb{N}_0$, $k=1,\ldots,2n+1$, i.e., $\sum_{n=0}^\infty\sum_{k=1}^{2n+1}n^4|\gamma_{n,k}|^2<\infty$. Here, $\delta_{n,m}$ denotes the Kronecker delta. 
\end{defi}
 
\begin{cor}[Existence of Induced Magnetizations]\label{thm:indmag2}
Let $m\in {\rm h}^2(\Omega)=\{f\in {\rm l}^2(\Omega):\sum_{i=1}^3\sum_{n=0}^\infty\sum_{k=1}^{2n+1}n^4|\langle f,\tilde{y}^{(i)}_{n,k}\rangle_{{\rm l}^2(\Omega)}|^2<\infty\}$ be a given magnetization. Then, for every admissible vector field $v\in{\rm l}^\infty(\Omega)$, there exists a $Q\in\LL^2(\Omega)$ such that the induced magnetization $\mathsf{m}=Qv$ is equivalent from outside to $m$. Analogously, for every admissible vector field $v\in{\rm l}^\infty(\Omega)$, there exists a $Q\in\LL^2(\Omega)$ such that the induced magnetization $\mathsf{m}=Qv$ is equivalent from inside to $m$.
\end{cor}

\begin{prf}
 According to Theorems \ref{thm:sunique1} and \ref{thm:shhdecomp2}, $m$ and $\mathsf{m}$ are equivalent from outside if and only if
 \begin{align}\label{eqn:hrel1}
  \DD^{-1}[\mathsf{M}_1]+\mathsf{M}_2+\frac{1}{2}\DD^{-1}[\mathsf{M}_2]=2\tilde{M}_2,
 \end{align}
 where $\mathsf{M}_1$, $\mathsf{M}_2$ are the Helmholtz scalars of $\mathsf{m}$ and $\tilde{M}_2$ the Hardy-Hodge scalar of $m$. Furthermore, since $v$ is admissible, we get $Q(\xi)=\frac{\mathsf{M}_1(\xi)}{\xi\cdot v(\xi)}$, $\xi\in\Omega$, and, therefore,
 \begin{align}
  \langle \mathsf{M}_2,Y_{n,k}\rangle_{\LL^2(\Omega)}&=-\frac{1}{n(n+1)}\langle \mathsf{M}_2,\Delta^*Y_{n,k}\rangle_{\LL^2(\Omega)}=\frac{1}{n(n+1)}\langle \nabla^*\mathsf{M}_2,\nabla^*Y_{n,k}\rangle_{{\rm l}^2(\Omega)}\label{eqn:hrel2}
  \\&=\frac{1}{n(n+1)}\langle \mathsf{m},\nabla^*Y_{n,k}\rangle_{{\rm l}^2(\Omega)}\nonumber
  \\&=\frac{1}{n(n+1)}\int_\Omega \mathsf{M}_1(\eta)\frac{v(\eta)}{\eta\cdot v(\eta)}\cdot\nabla^*Y_{n,k}(\eta)d\omega(\eta)\nonumber
  \\&=\frac{1}{n(n+1)}\sum_{m=0}^\infty\sum_{l=1}^{2m+1}\langle \mathsf{M}_1,Y_{m,l}\rangle_{\LL^2(\Omega)}\int_\Omega Y_{m,l}(\eta)\frac{v(\eta)}{\eta\cdot v(\eta)}\cdot\nabla^*Y_{n,k}(\eta)d\omega(\eta)\nonumber
  \\&=\frac{1}{n(n+1)}\sum_{m=0}^\infty\sum_{l=1}^{2m+1}\langle \mathsf{M}_1,Y_{m,l}\rangle_{\LL^2(\Omega)}\,v_{n,k,m,l}.\nonumber
 \end{align}
 Observing that the operator $\DD^{-1}$ acts via
 \begin{align}
  \DD^{-1}[\mathsf{M}_1]=\sum_{n=0}^\infty\sum_{k=1}^{2n+1}\frac{1}{n+\frac{1}{2}}\langle \mathsf{M}_1,Y_{n,k}\rangle_{\LL^2(\Omega)}Y_{n,k}
 \end{align}
 and using \eqref{eqn:hrel2} in \eqref{eqn:hrel1} leads to the following infinite dimensional system of linear equations for the Fourier coefficients $\langle \mathsf{M}_1,Y_{n,k}\rangle_{\LL^2(\Omega)}$ of $\mathsf{M}_1$:
\begin{align}
 \sum_{m=0}^\infty\sum_{l=1}^{2m+1}\left(\frac{1}{n(n+\frac{1}{2})}v_{n,k,m,l}+\frac{1}{n+\frac{1}{2}}\delta_{n,m}\delta_{k,l}\right)\langle \mathsf{M}_1,Y_{m,l}\rangle_{\LL^2(\Omega)}=&2\langle \tilde{M}_2,Y_{n,k}\rangle_{\LL^2(\Omega)},
 \\[-2ex]&\qquad n\in\mathbb{N}_0,k=1,\ldots,2n+1.\nonumber
\end{align}
The admissibility conditions on $v$ guarantee that this problem is solvable and that the obtained $Q(\xi)=\frac{\mathsf{M}_1(\xi)}{\xi\cdot v(\xi)}$, $\xi\in\Omega$, lies in $\LL^2(\Omega)$. The Fourier coefficients $\langle \mathsf{M}_2,Y_{n,k}\rangle_{\LL^2(\Omega)}$ of $\mathsf{M}_2$ can be obtained via \eqref{eqn:hrel2}. Here, $m\in{\rm h}^2(\Omega)$ guarantees a sufficient decay of the Fourier coefficients  $\langle \mathsf{M}_1,Y_{n,k}\rangle_{\LL^2(\Omega)}$ such that the interchange of the series and integration in \eqref{eqn:hrel2} is allowed. The resulting induced magnetization $\mathsf{m}=Qv$ is then equivalent from outside to $m$.
\end{prf}

\begin{rem}
If we choose $v(\xi)=\xi$, $\xi\in\Omega$ (which served as a counterexample in Remark \ref{rem:locexist} for the case supp$(\mathsf{m})\subset\Gamma$), the conditions on $v$ of Definition \ref{def:adm} are clearly satisfied since $v_{n,k,m,l}=0$. In other words, for any magnetization $m\in{\rm l}^2(\Omega)$ there exists an induced magnetization $\mathsf{m}=Qv$ that is equivalent from outside. For general vector fields $v$ it is more difficult to check the conditions from Definition \ref{def:adm}. Here, the reader is, e.g., referred to \cite{shivakumar09} and references therein. 
\end{rem}

\section{Numerical Illustrations}\label{sec:num}

\begin{figure}
\begin{center}
\includegraphics[scale=0.33]{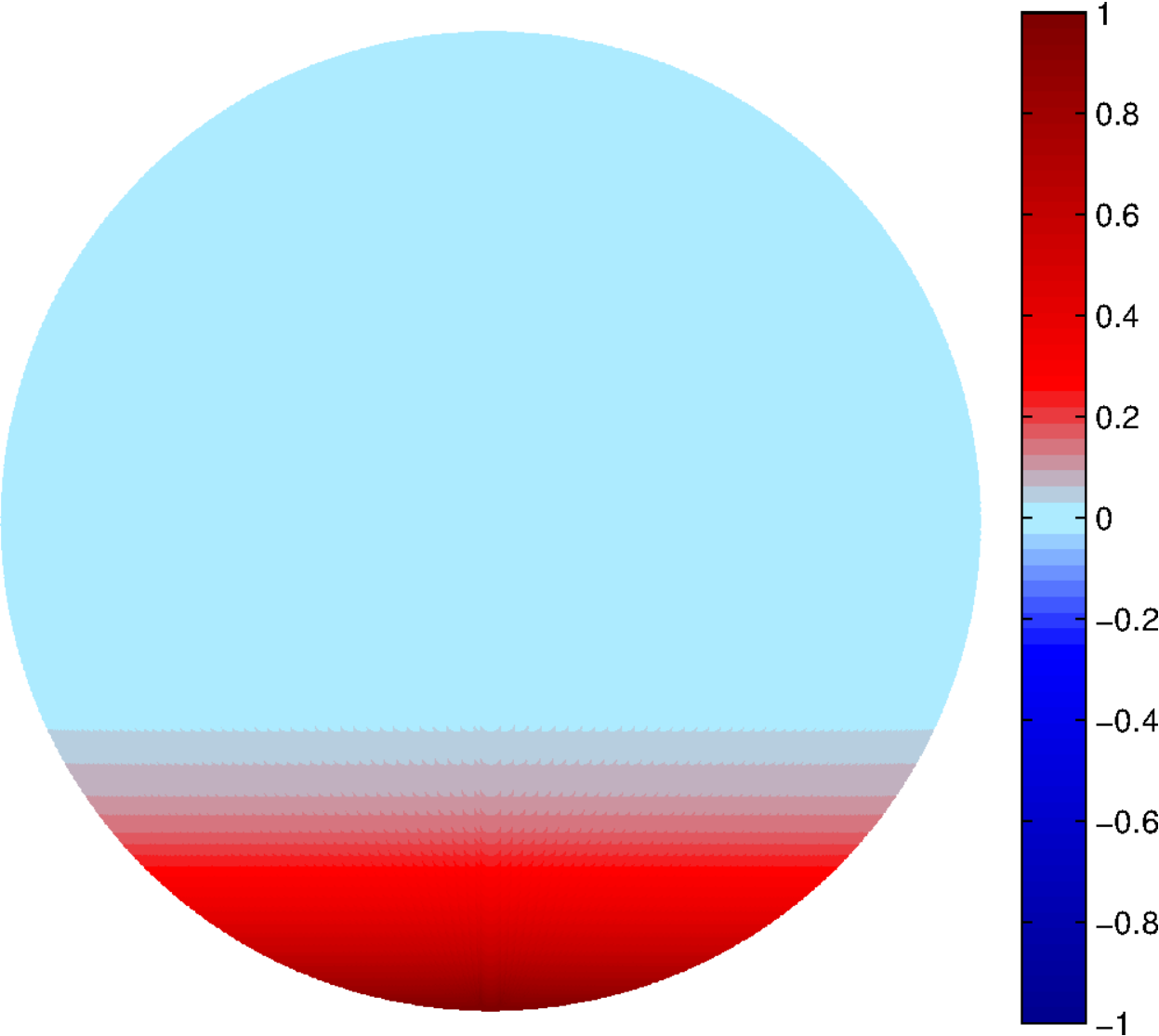}\qquad\qquad\qquad\includegraphics[scale=0.21]{sunidi1_2-eps-converted-to.pdf}
\end{center}
\caption{The true susceptibility $Q$ (left) and the inducing vector field $v$ (right).}
\label{fig:truem}
\begin{center}
\includegraphics[scale=0.37]{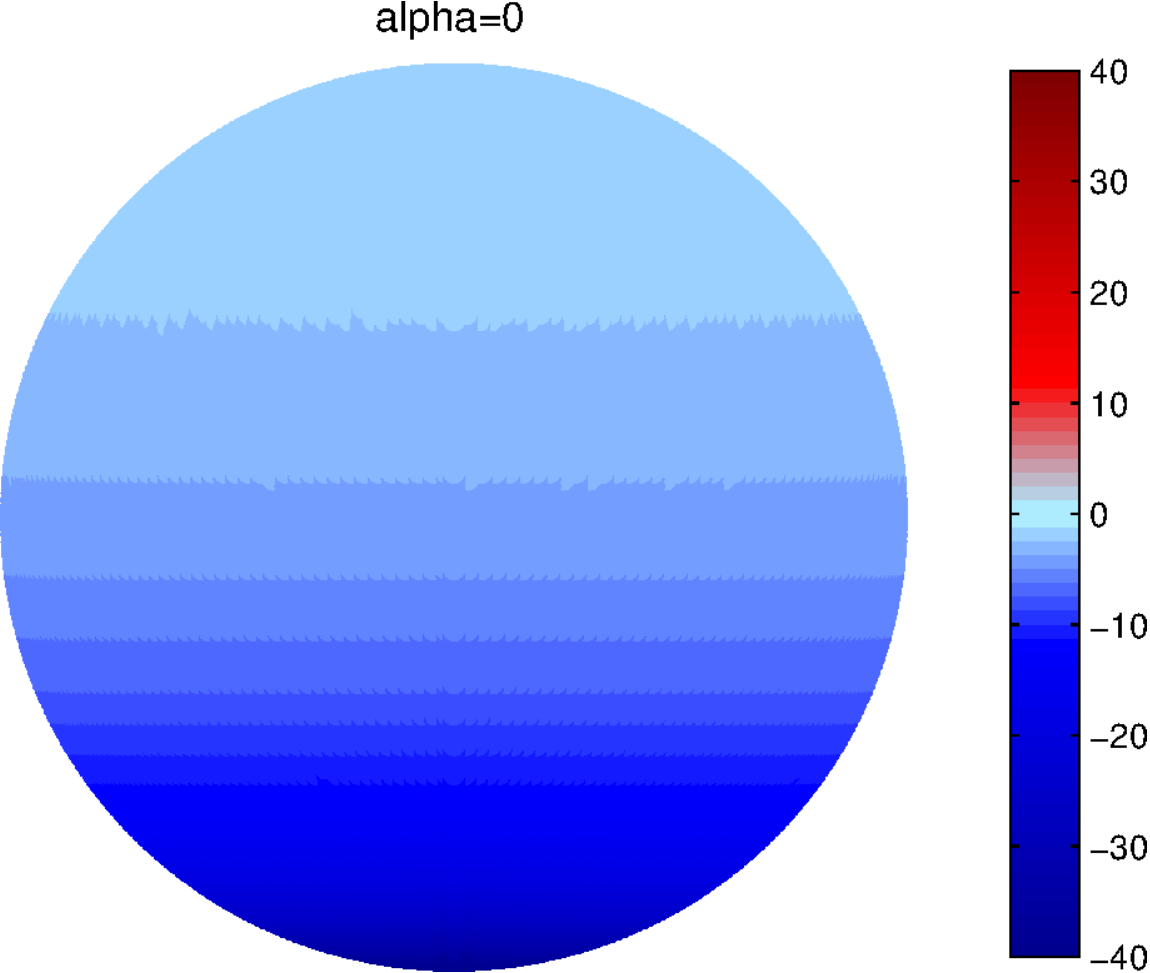}\qquad\includegraphics[scale=0.3]{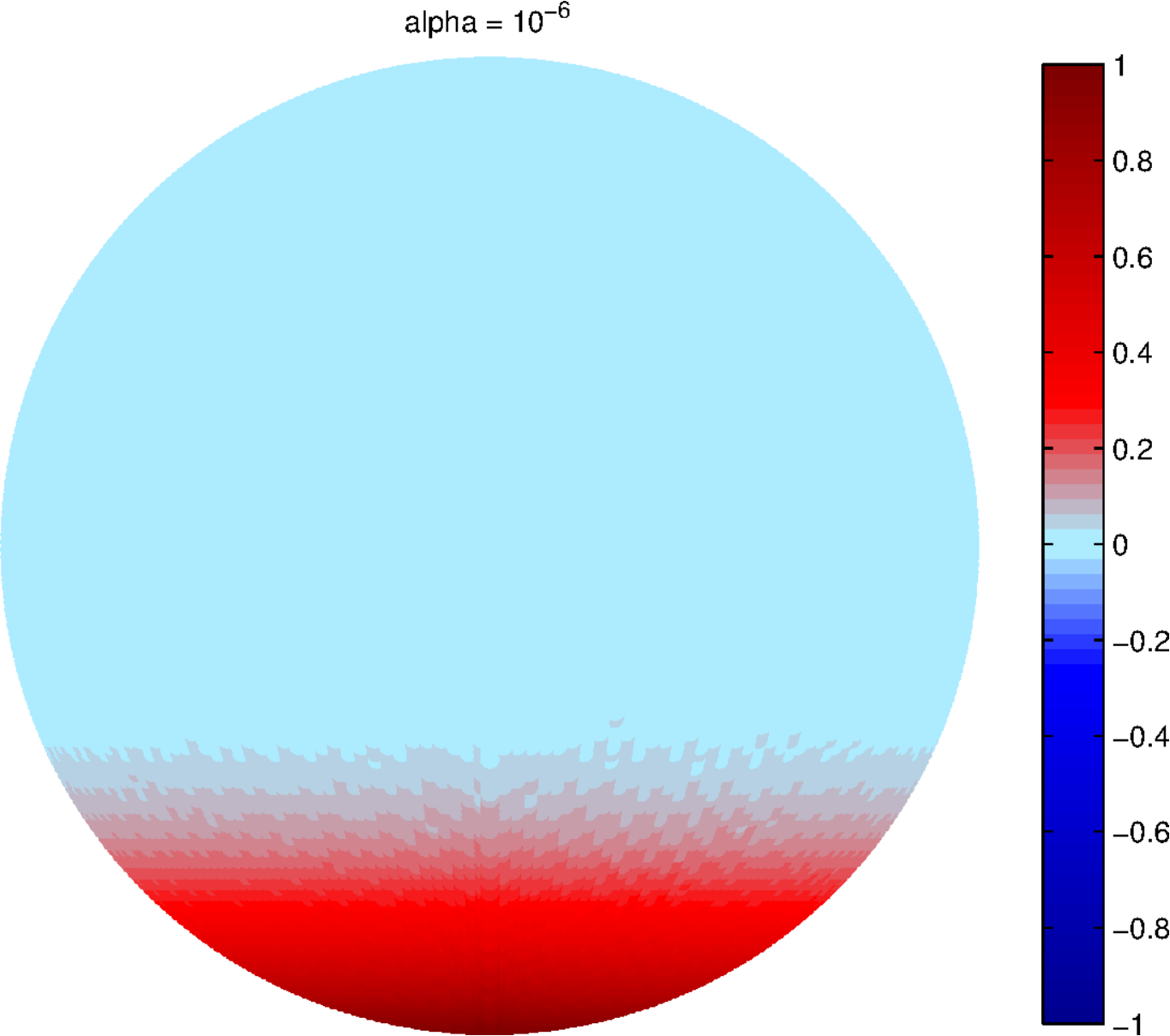}\qquad\includegraphics[scale=0.3]{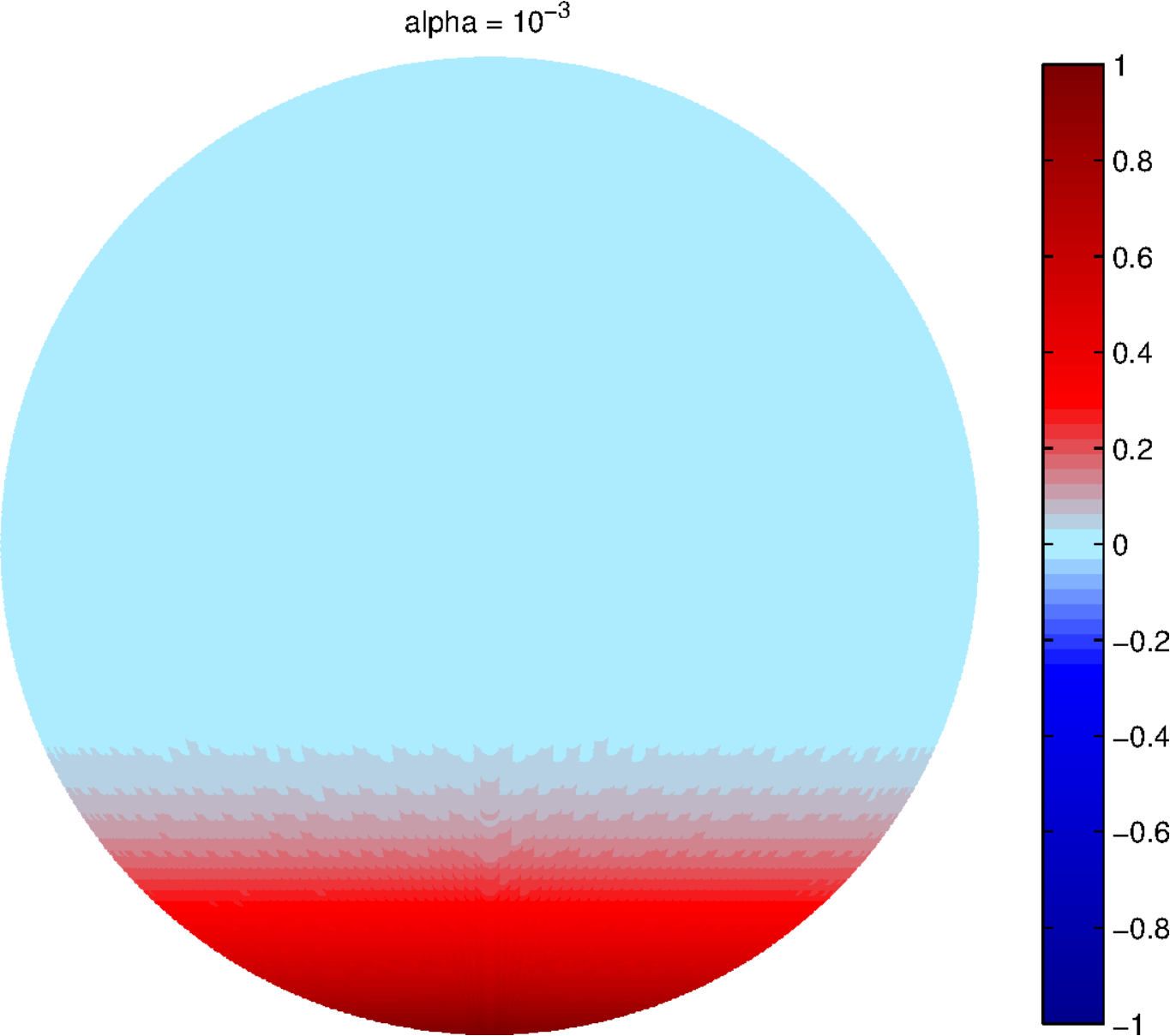}
\end{center}
\caption{The reconstructed susceptibilities $\bar{Q}$ for different values $\alpha=0,10^{-6},10^{-3}$ (the images represent a frontal view of the sphere on which $\bar{Q}$ is plotted).}
\label{fig:reconstm}
\end{figure}

In this section, we want to illustrate that the theoretical result from Corollary \ref{thm:indmag} has an actual influence on the numerical reconstruction of induced magnetizations $m=Qv$. We assume to know the magnetic potential $V$ on $\Omega_R$ (we choose $R=1.1$ in this example) that is generated by the induced magnetization $m=Qv$, where 
\begin{align}
 Q(\xi)&=\left\{\begin{array}{ll}(\xi\cdot\zeta)^4,&\xi\cdot\zeta\leq0,
 \\0,&\textnormal{else},\end{array}\right.
 \\v(\xi)&=\xi+\zeta-(\zeta\cdot\xi)\xi,\label{eqn:vind}
\end{align}
and $\zeta=(0,0,1)$ is fixed, i.e., the magnetization $m$ has compact support in the lower hemisphere $\Gamma=\{\xi\in\Omega:\xi\cdot\zeta\leq0\}$ (see Figure \ref{fig:truem} for an illustration). 

In order to approximate the true (but unknown) magnetization $m=Qv$ by some $\bar{m}=\bar{Q}v$, we denote by $V[\bar{Q}]$ the magnetic potential that is generated by $\bar{m}$ and minimize the functional
\begin{align}\label{eqn:minfunc}
 \mathcal{F}[\bar{Q}]=\left\|V[\bar{Q}]-V\right\|_{\LL^2(\Omega_R)}^2+\alpha\|\bar{Q}v\|_{{\rm l}^2(\Omega\setminus\Gamma)}^2.
\end{align}
The first term in \eqref{eqn:minfunc} simply represents a data misfit that measures the deviation of $V[\bar{Q}]$ from the known magnetic potential $V$. The second term in \eqref{eqn:minfunc} penalizes magnetizations $\bar{m}=\bar{Q}v$ that have contributions outside $\Gamma$, i.e., magnetizations that do not satisfy supp$(\bar{m})\subset\Gamma$. For the numerical minimization of $\mathcal{F}[\bar{Q}]$, we expand $\bar{Q}$ in terms of Abel-Poisson kernels (cf. \cite{freeden98}):
\begin{align}
 \bar{Q}(\xi)&=\sum_{n=1}^N\gamma_n K(\xi\cdot \xi_n),
 \\K(\xi\cdot\xi_n)&=\frac{1-h^2}{(1+h^2-2h(\xi\cdot\xi_n))^{\frac{3}{2}}},
\end{align}
where $h\in(0,1)$ is a fixed parameter (influencing the localization of $K$) and $\xi_n\in\Omega$, $n=1,\ldots,N$, are some predefined points indicating different centers of the kernel $K$ (in our case, we choose $N=10235$ and points $\xi_n$ that are uniformly distributed on $\Omega$). Under these conditions, the minimization of $\mathcal{F}[\bar{Q}]$ reduces to solving the set of linear equations
\begin{align}
 \mathbf{M}\gamma=g,
\end{align}
where
\begin{align}
 \mathbf{M}&=(M_{n,m})_{n,m=1,\ldots,N}
 \\&=\left(\int_{\Omega_R}V_n(x)V_m(x)d\omega(x)+\alpha \int_{\Omega\setminus\Gamma}|v(\eta)|^2K(\eta\cdot\xi_n)K(\eta\cdot\xi_m)d\omega(\eta)\right)_{n,m=1,\ldots,N},\nonumber
 \\\gamma&=(\gamma_m)_{m=1,\ldots,N},
 \\g&=(g_n)_{n=1,\ldots,N}=\left(\int_{\Omega_R}V_n(x)V(x)d\omega(x)\right)_{n=1,\ldots,N},
  \end{align}
and
\begin{align}
V_n(x)&=\frac{1}{4\pi}\int_\Omega K(\eta\cdot\xi_n)v(\eta)\cdot\frac{x-\eta}{|x-\eta|^3}d\omega(\eta).\label{eqn:V_n}
\end{align}
The quadrature rules from \cite{driscoll94,hesse12} are used for the numerical evaluation of the occurring integrals. The reconstructed susceptibilities $\bar{Q}$ for different choices of $\alpha$ are shown in Figure \ref{fig:reconstm} (the parameter $h$ of the Abel-Poisson kernels is set to $h=0.9$). One can see that for parameters $\alpha>0$, we obtain a fairly good reconstruction of the underlying true susceptibility $Q$, while for $\alpha=0$ (i.e., no penalization is taken into account for magnetizations $\bar{m}$ that violate supp$(\bar{m})\subset\Gamma$), we obtain an entirely different susceptibility $\bar{Q}$. Latter generates the same magnetic potential on $\Omega_R$ as $Q$ but it does not satisfy supp$(\bar{m})\subset\Gamma$. 

\begin{figure}
\begin{center}
\includegraphics[scale=0.35]{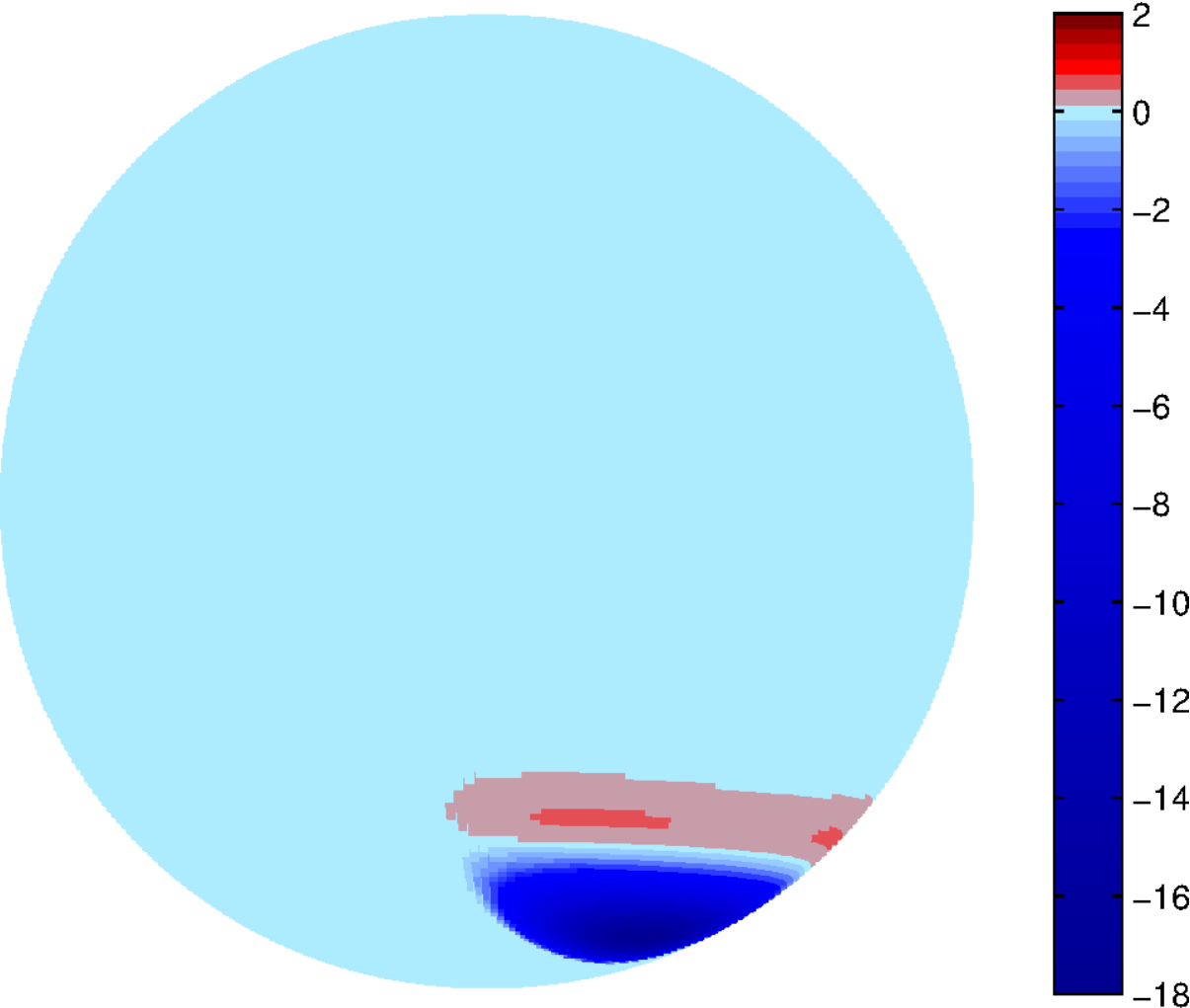}\qquad\qquad\qquad\includegraphics[scale=0.28]{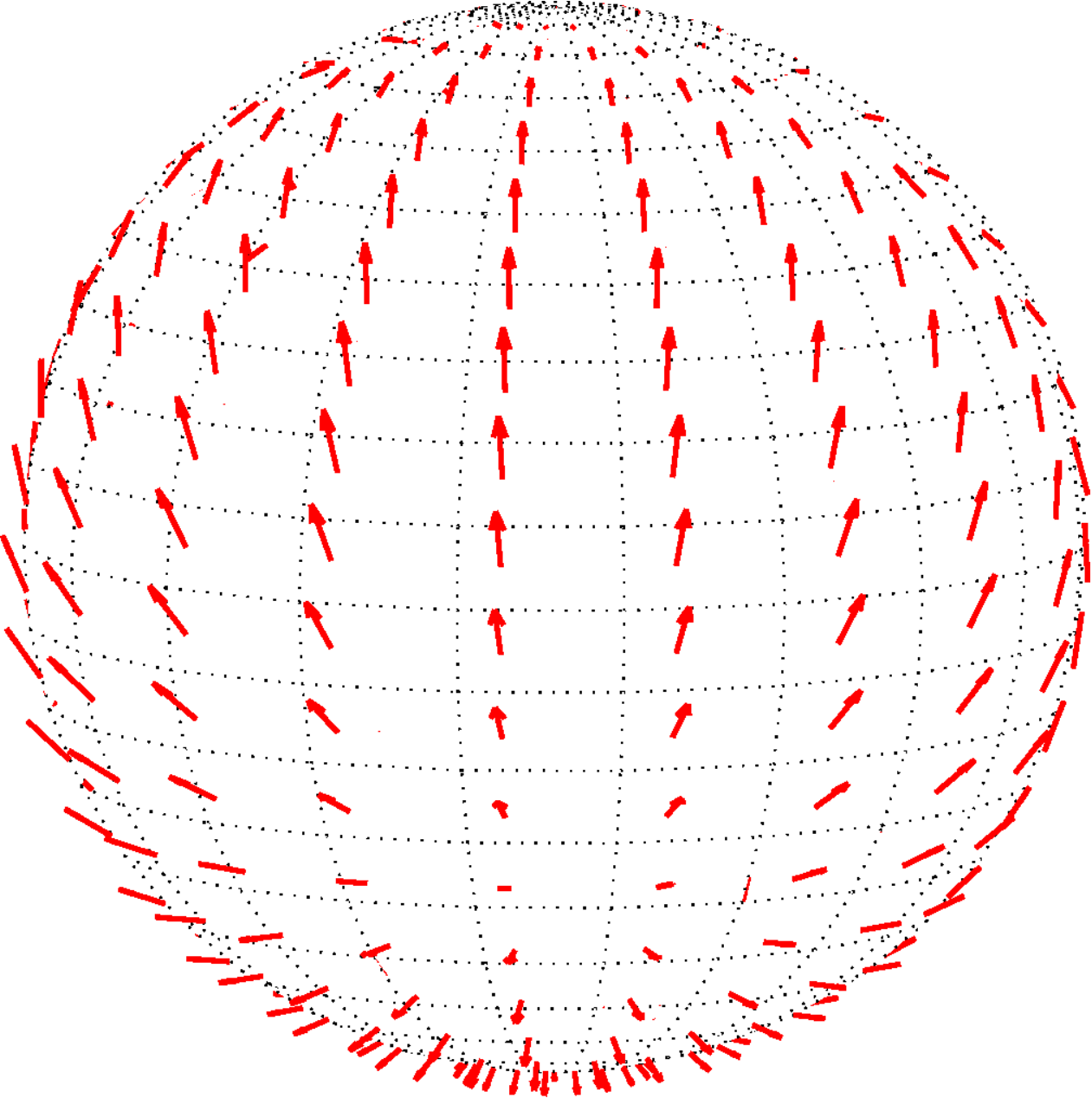}
\end{center}
\caption{The true susceptibility $Q$ (left) and the inducing vector field $v$ (right).}
\label{fig:truem2}
\begin{center}
\includegraphics[scale=0.36]{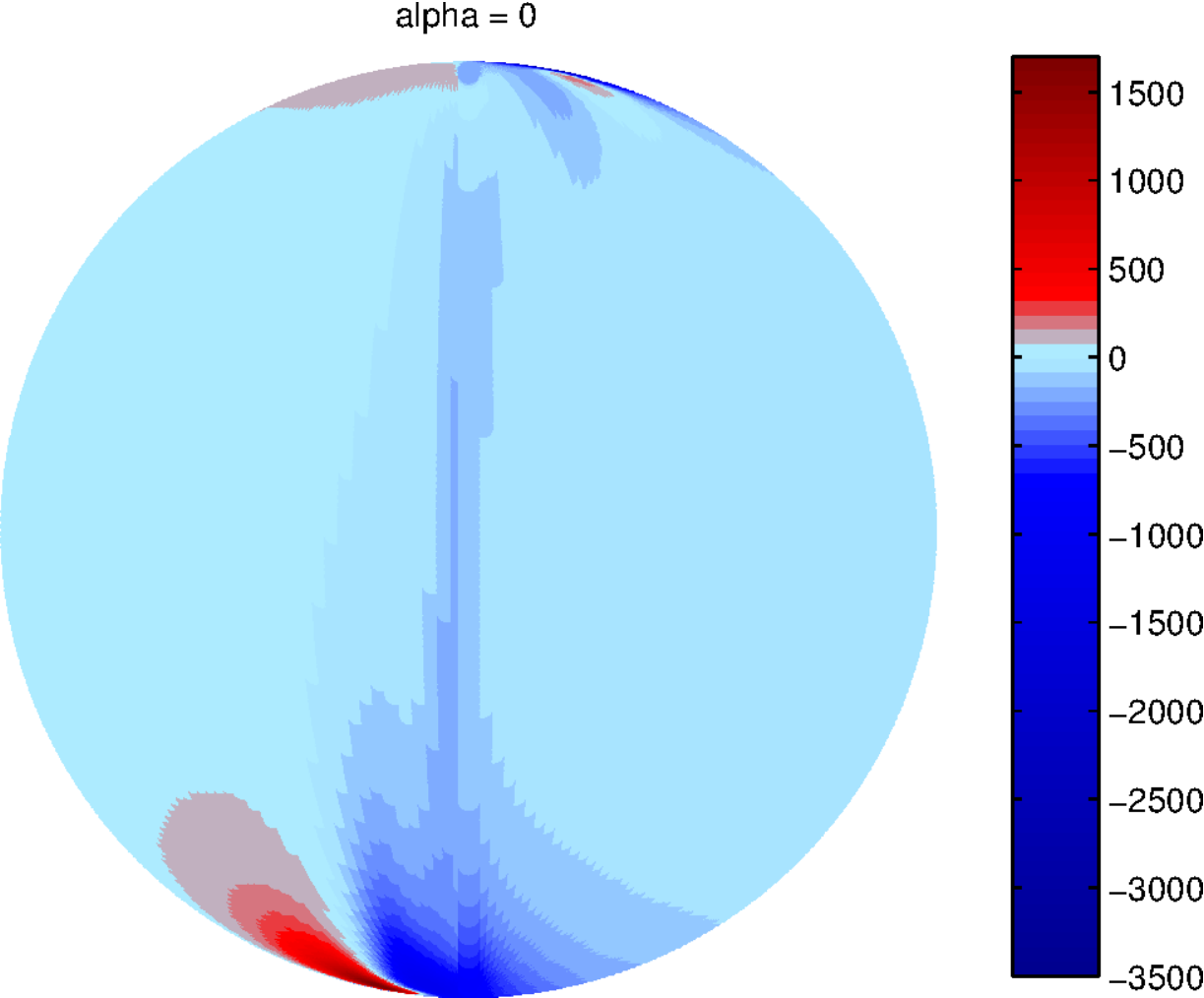}\qquad\includegraphics[scale=0.37]{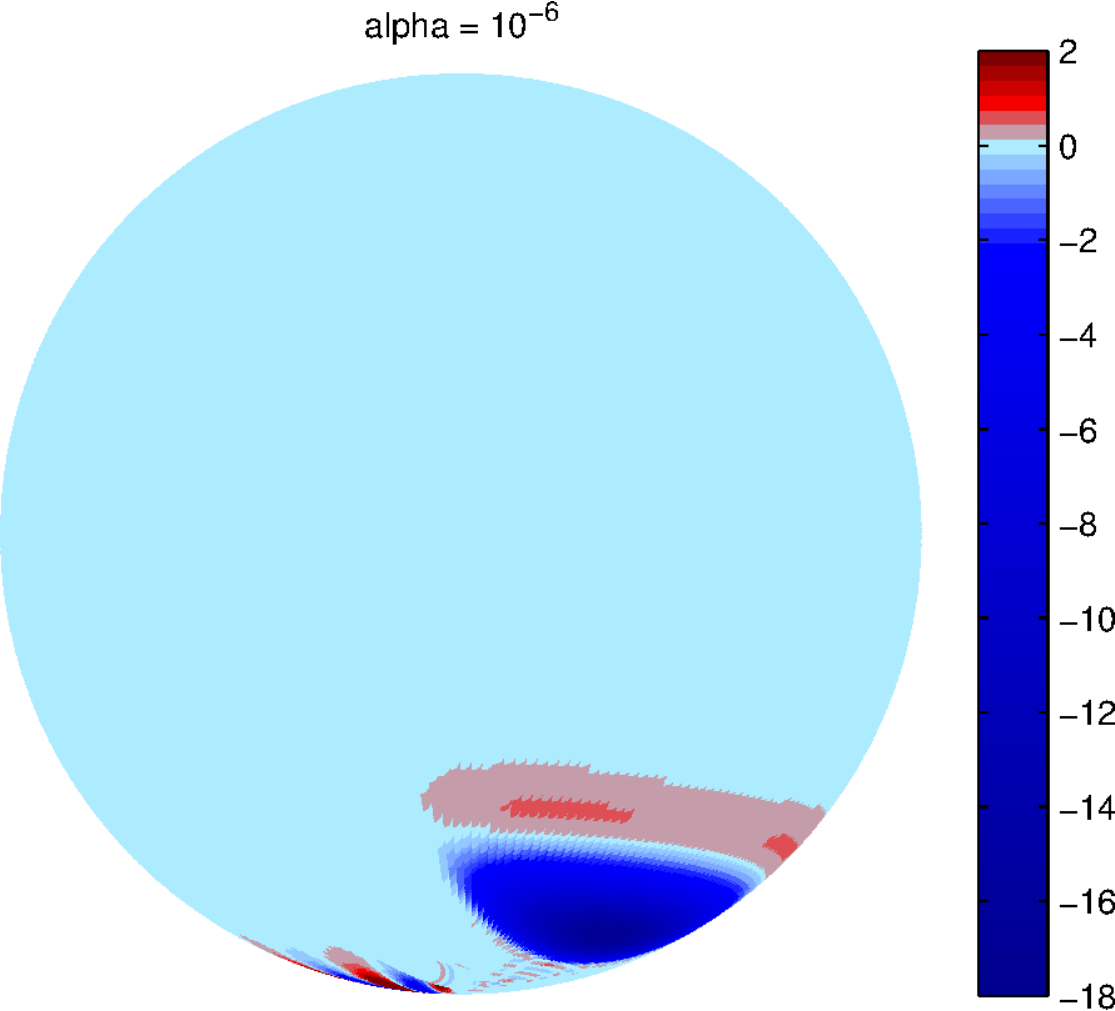}\qquad\includegraphics[scale=0.4]{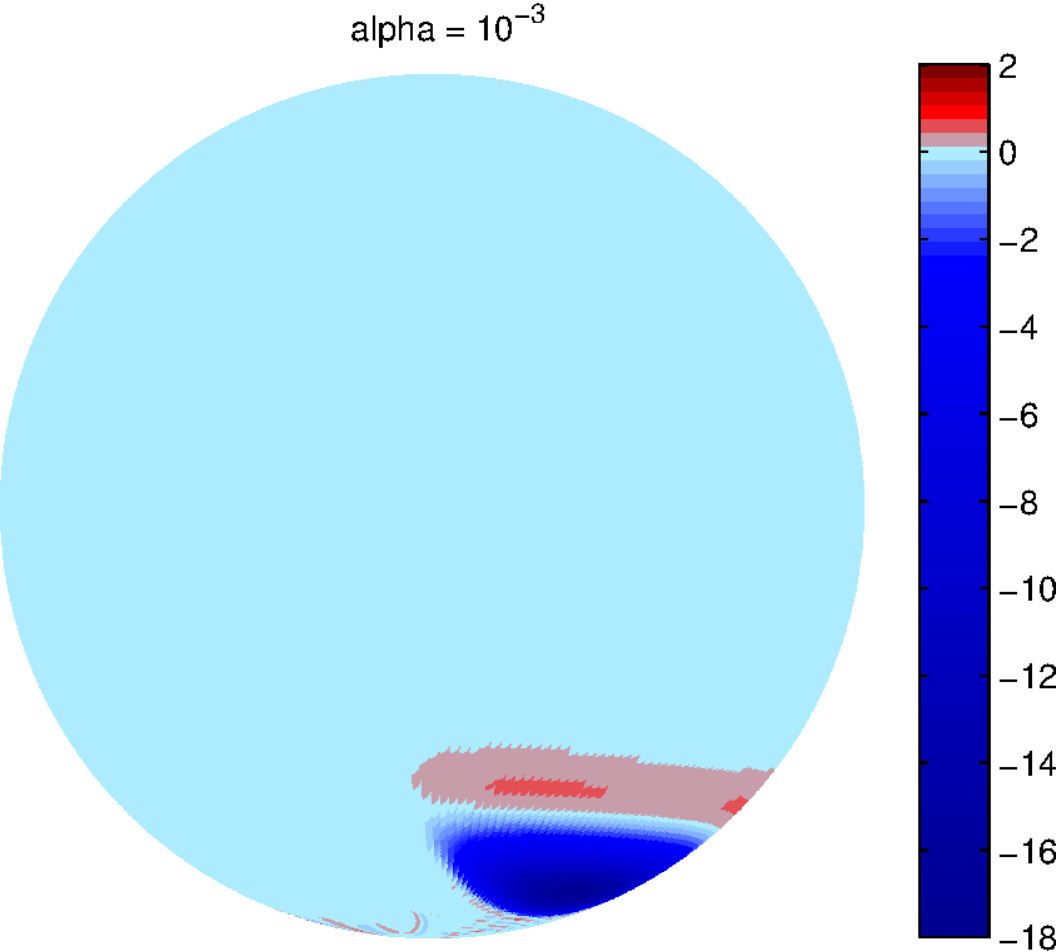}
\end{center}
\caption{The reconstructed susceptibilities $\bar{Q}$ for different values $\alpha=0,10^{-6},10^{-3}$ (the images represent a frontal view of the sphere on which $\bar{Q}$ is plotted).}
\label{fig:reconstm2}
\end{figure}

In a second test example, we choose a slightly more complicated magnetization $m=Qv$:
\begin{align}
Q(\xi)&=\left\{\begin{array}{ll}1000\left(\frac{1}{2}+\xi\cdot\zeta\right)^4\cos(2\pi\xi\cdot\zeta)\sin(2\pi \xi\cdot\bar{\zeta}),&\xi\cdot\zeta\leq-\frac{1}{2}\textnormal{ and }\xi\cdot\bar{\zeta}\geq\frac{1}{2},
\\0,&\textnormal{else},\end{array}\right.
\\v(\xi)&=(\bar{\zeta}\cdot\xi)\xi+\zeta-(\zeta\cdot\xi)\xi,
\end{align}
where $\zeta=(0,0,1)$ and $\bar{\zeta}=(0,\frac{\sqrt{15}}{4},-\frac{1}{4})$ are fixed. The magnetization $m$ is again supported in the lower hemisphere $\Gamma=\{\xi\in\Omega:\xi\cdot\zeta\leq0\}$, but the actual support is only a subset of $\Gamma$ (see Figure \ref{fig:truem2} for an illustration). $V$ denotes the magnetic potential on $\Omega_R$ produced by $m$. In order to approximate the (unknown) magnetization $m=Qv$ by a magnetization $\bar{m}=\bar{Q}v$ from knowledge of $V$, we again minimize $\mathcal{F}[\bar{Q}]$ as in \eqref{eqn:minfunc}--\eqref{eqn:V_n}. The parameter $h$ of the involved Abel-Poisson kernels is chosen to be $h=0.95$ (in order to supply better localized kernels $K$). The results are shown in Figure \ref{fig:reconstm2}. Once more, for $\alpha=0$, we obtain a magnetization $\bar{m}$ that is equivalent to $m$ from outside but does not satisfy supp$(\bar{m})\subset\Gamma$. Choosing $\alpha>0$ leads to a reconstruction of the desired magnetization. However, we see (in particular in the center image of Figure \ref{fig:reconstm2}) that there are some undesired artefacts around the South Pole. Considering that these artefacts are reduced for an increased $\alpha$, they might be due to the fact that the operator that maps the magnetization $m$ on $\Omega$ to its magnetic potential $V$ on $\Omega_R$ has an unbounded inverse (the classical ill-posedness that typically occurs for potential field problems). In this paper, our focus is on ill-posedness in the sense of non-uniqueness. Some more sophisticated regularization methods that also address the ill-posedness due to unboundedness of the inverse for the planar magnetization problem can be found, e.g., in \cite{lima13}.

Last, it should be noted that in the second example, the lower hemisphere $\Gamma$ that we chose for the numerical reconstruction is significantly larger than the actual support of $m$. This shows that the numerical scheme also works if we do not know the exact support of $m$.

\section{Conclusion}
We proved that for induced spherical magnetizations (where the inducing vector field is known) the additional assumption of compact support in some region $\Gamma\subset\Omega$ yields uniqueness for $m$. The numerical examples indicate that including this additional condition in the reconstruction procedure guarantees picking the 'correct' magnetization out of those that could generate the measured magnetic potential.


\end{document}